\RequirePackage{fix-cm}
\documentclass[smallextended]{svjour3}       
\usepackage[utf8]{inputenc}

\usepackage{tabularx,ragged2e,booktabs,caption}
\newcolumntype{C}[1]{>{\Centering}m{#1}}

\smartqed  
\usepackage{xspace}
\usepackage{amsfonts}
\usepackage{amsmath}
\usepackage{enumerate}
\usepackage{framed}
\usepackage{graphicx}
\usepackage{lscape}
\usepackage{bm}
\usepackage{color}
\usepackage{multicol}
\usepackage[utf8]{inputenc}

\usepackage{tabularx,ragged2e,booktabs,caption}
\newcolumntype{C}[1]{>{\Centering}m{#1}}

\usepackage{amsmath}
\usepackage{amsfonts}
\usepackage{amssymb}
\usepackage{tikz}
\usepackage{subfig}
 \usepackage{float}
\usepackage{slashbox}
\usepackage{graphicx}
\usepackage[english]{babel}
\usepackage{geometry}                
\usepackage[parfill]{parskip}    
\usepackage{epstopdf}
\usepackage{psfrag}
\newcommand{\figref}[1]{{Figure~\ref{#1}}}
\usepackage[sort]{natbib}
\newcommand{\tabref}[1]{{Table~\ref{#1}}}

\renewcommand{\eqref}[1]{{(\ref{#1})}}

 \newcommand{\M}{\mathbf{M}}
 \newcommand{\C}{\mathcal{C}}
 \newcommand{\U}{\mathcal{U}}
\newcommand{\R}{\mathcal{R}}
\newcommand{\V}{\mathcal{V}}
\newcommand{\F}{\mathcal{F}}
\newcommand{\E}{\mathcal{E}}

\usepackage{graphicx}
\begin{document}

\title{ A Fitted Multi-Point Flux Approximation Method  for  Pricing two options}
%


\titlerunning{Fitted Multi-Point Flux Approximation method  for pricing options}        

\author{ Rock Stephane Koffi,
Antoine Tambue
}

\authorrunning{ R. S. Koffi, A. Tambue} 

\institute{A. Tambue (Corresponding author) \at
             Western Norway University of Applied Sciences,  Inndalsveien 28, 5063  Bergen, Norway,\\
            The African Institute for Mathematical Sciences(AIMS), 6-8 Melrose Road, Muizenberg 7945, South Africa\\
            Center for Research in Computational and Applied Mechanics (CERECAM), and Department of Mathematics and Applied Mathematics, University of Cape Town, 7701 Rondebosch, South Africa.\\
             Tel.: +47 55 58 70 06, \email{antonio@aims.ac.za, antoine.tambue@hvl.no, tambuea@gmail.com} \\
   \and
   R.S. Koffi \at
            The African Institute for Mathematical Sciences(AIMS) , 6-8 Melrose Road, Muizenberg 7945, South Africa\\
            Department of Mathematics and Applied Mathematics, University of Cape Town, 7701 Rondebosch, South Africa\\
            \email{rock@aims.ac.za}          
}
\date{Received: date / Accepted: date}

\maketitle

\begin{abstract}
In this paper, we develop novel numerical methods based on the Multi-Point Flux Approximation (MPFA) method to solve the degenerated partial differential equation (PDE) arising 
from pricing two-assets options. The standard   MPFA   is used  as our  first method  and   is    coupled  with a fitted finite volume  in our second  method to handle the degeneracy of the PDE  and the corresponding scheme is called fitted MPFA method.
The convection part is discretized using the upwinding methods (first and second order)  that we have derived on non uniform grids.
The  time discretization  is performed with $\theta$- Euler methods.
Numerical simulations show that our new schemes  can be  more accurate  than the current fitted finite volume method  proposed in the literature.

\keywords{Finite volume methods, Multi-Point Flux Approximation,  Degenerated  PDEs, Options pricing, Multi-asset options }
  
\end{abstract}

\section{Introduction}
\label{intro}
 Pricing multi-assets options  is of great interest in the financial industry (see  \cite{persson2007pricing}).
Multi-asset options are options based on more than one underlying. There are several kinds of multi-assets options, few  of  them are  exchange options, rainbow options, 
baskets options, best or worst options, quotient options, foreign exchange options, quanto options, spread options, dual-strike options and out-performance options. 
Pricing these options lead to the resolution of the following second order degenerated Black-Scholes Partial Differential Equations (PDE)(see \cite{persson2007pricing})

\begin{eqnarray}
\label{multi}
\frac{\partial U}{\partial\tau}=\frac{1}{2}\sum_{i,j=1}^n \sigma_i\sigma_j\rho_{ij}S_iS_j\frac{\partial^2U}{\partial S_i \partial S_j}+r\sum_{i=1}^nS_i\frac{\partial U}{\partial S_i}-rU
\end{eqnarray}
where $r$ is the risk free interest, $U$ is the option  value at time $\tau$, $\tau=T-t$ with $t$ and $T$ respectively the instantaneous and  maturity time,
$\,S_i~$ represents the asset $i$ price, $\sigma_i$  represents the volatility of asset $i$,  $\rho_{ij}$  represents the correlation between the assets $i$ and $j$, where $i, j=1,...,n$.
The main difference between multi-assets options is their payoff functions which represent the initial condition of the corresponding backward PDE.
The spatial domain of the PDE is infinite, but for its numerical resolution, a truncation is required (see  \cite{duffy2013finite},Chapter 3). It has been observed that when the stock price $S$ approaches 
the region near to zero, the Black Scholes PDE is degenerated (see \cite{duffy2013finite}, chapter 30.3). Moreover, the initial condition of the PDE has a  discontinuity in its first derivative 
when the stock price is equal to the strike $K$. This discontinuity has an adverse impact on the accuracy when the finite difference method is used (see \cite{wilmott2005best}, chapter 26). Therefore, 
for the spatial  discretization of the  PDE, it is suitable to use non-uniform grids with more points in the region around $S=0$ and $S=K$ in order to handle the degeneracy and the discontinuity.
To overcome the above challenges, many methods have been proposed in the literature. Thereby, \cite{wang2004novel} proposed  a  fitted  finite volume method  for one dimensional Black Scholes PDE 
 and the rigorous convergence proof is provided by \cite{angermann2007convergence}. Besides, \cite{huang2006fitted} adapted  the fitted finite volume discretization method for the two-dimensional Black-Scholes PDE   and  its rigorous convergence proof
is analysed by \cite{huang2009convergence}. Although these two fitted finite volume methods are stable,  they are only order 1 with respect to asset price variables.           

In this paper, we present two novel discretization methods for the two-dimensional Black Scholes PDE based on a special kind of  finite volume method,
the so-called  Multi-Point Flux Approximation (MPFA) method. This method  was introduced by \cite{aavatsmark2002introduction} and has been  used in fluid dynamics  
for flow and transport equations (see \cite{sandve2012efficient} and references therein). Actually, the MPFA was designed to give a correct discretization of the flow equation for general  
grids including fractures (see \cite{aavatsmark2002introduction,sandve2012efficient}). 
The MPFA  method is essentially based on the approximation of a linear function gradient over a triangle, the calculation and the continuity of  
flux through edges of this triangle. The convergence of  MPFA method is usually second order in space domain  on rough grids (see \cite{aavatsmark2007multipoint,stephansen2012convergence}).
Our first numerical method  here is  the  standard MPFA ,  which  is  fully  used to approximate the second order operator.  
To the best  of our knowledge,  this method was not yet used to solve degenerated Black Scholes PDE  in finance.
To build our new fitted  MPFA method, we couple the standard MPFA with 
the upwind methods (first and second order) to approximate two dimensional options pricing.
Besides, the fitted finite volume proposed by \cite{wang2004novel} is used to handle the degeneracy of the PDE in the region where the stocks price approach zero (degeneracy region). In the region, where the PDE
in not degenerated, we apply the MPFA method. The  novel numerical technique from  this combination is called fitted MPFA method and will obviously
improve the accuracy of the current fitted finite volume in the literature, since more approximations involving are second order in space.  Naturally, 
these  two methods  are  applicable to  other types of multi-asset options and also to financial models such as  \cite{heston1993closed} model and  \cite{bates1996jumps} model  on non-uniform grids. Another advantage of our novel fitted MPFA is that it can easily be adapted to more structured commercial or open-source softwares as the standard MPFA (see \cite{lie2012open}). \\
The rest of the paper is organized as follows. In section 2, we start by introducing the Black Scholes model for option with 2 stocks and the 
corresponding partial differential equation. Afterwards, we set the frame of  the numerical domain of study suitable for the finite volume method application.
Section 3 is devoted to the spatial discretization of the PDE. We describe the Multi-Point Flux Approximation method for the discretization of the diffusion term of the PDE.
The upwind methods (first and second order) are used for the the convection term discretization. We end the section 3 with the fitted MPFA which is a combination of 
a fitted finite volume method and the MPFA method. The time discretization is performed using the $\theta-$Euler methods in the section 4. In section 5, we perform numerical 
experiments. Those numerical simulations show that the two proposed  schemes (the standard MPFA  method and fitted MPFA method ) can be  more accurate  than the current fitted finite volume method  proposed in the literature.
General conclusion is given in section 6.

%
\section{Formulation of the problem}

\subsection{Black-Scholes model with 2 underlying assets}

An option with two underlying assets modeled by the Black Scholes equation is formulated as follows

\begin{align}
\left\lbrace \begin{array}{l}
dx(t)~~~~=~\mu_1xdt+\sigma_1 xdW_1\\
\\
dy(t)~~~~=~\mu_2ydt+\sigma_2 ydW_2\\
\\
dW_1(t)dW_2(t)  =~\rho dt
\end{array}\right.
\end{align}

where $\mu_i,\sigma_i, W_i$ are respectively the drift, the volatility and the Wiener process governing the stocks $x,y$ and $\rho$ is 
the correlation coefficient between the two Wiener processes.
By applying  the Ito's formula and using the standard  arbitrage argument, it is well known ( see \cite{hull2003options,kwok2008mathematical,wilmott1993option} ) that  the value of the option  $U$ follows 
the following two-dimensional Black-Scholes Partial differential equation 
on the domain $D=[0,+\infty) \times[0,+\infty) \times[0,T]$
\begin{equation}
\label{twoop}
 \frac{\partial U}{\partial \tau}= \frac{1}{2} \sigma^2_1 x^2 \frac{\partial^2 U}{\partial x^2} +\rho \sigma_1 \sigma_2 xy\frac{\partial^2 U}{\partial x \partial y} + \frac{1}{2} \sigma_2^2 y^2 \frac{\partial^2 U}{\partial y^2}+rx\frac{\partial U}{\partial x} +ry\frac{\partial U}{\partial y}-rU
\end{equation} 

where $\tau=T-t$, $T$ is the maturity time, $t$ the current time and $r$ is the risk-free interest. 
 For European rainbow option price on maximum of two risky assets, the following  initial and boundary conditions are used

\begin{align}
\left\lbrace
\begin{array}{l}
U(x,y,0)=\max\left(\max(x,y)-K,0\right) \\
\\
U(0,y,\tau)=0\\
\\
U(x,0,\tau)=0\\
\end{array}\right.
\end{align}
with $K$ the strike price. But  to  compare our numerical solution with the existing fitted finite volume method, the exact solution  will be used at the boundary.  In order to apply the finite volume method, it is convenient to  re-write the Partial Differential Equation \eqref{twoop} in the following divergence form
\begin{eqnarray}
\label{conservation}
\frac{\partial U}{\partial \tau }= \nabla \cdot(\M\nabla U)+\nabla (f U)+\lambda U
\end{eqnarray}
where

\begin{eqnarray*}
& & \M=\frac{1}{2}\left(\begin{array}{lr}
\sigma_1^2 x^2 & \rho\sigma_1\sigma_2xy \\
   & \\
\rho\sigma_1\sigma_2xy & \sigma_2^2 y^2
\end{array}\right),
f=\left(\begin{array}{c}
(r-\sigma_1^2-\frac{1}{2}\rho\sigma_1\sigma_2)x \\ \\ (r-\sigma_2^2-\frac{1}{2}\rho\sigma_1\sigma_2)y
\end{array}\right)\\
& & \\
&  &\\
&   &~~~~~~~~~~~~~~~~~~~~~~~\lambda = -3r+\sigma_1^2+\sigma^2_2+\rho\sigma_1\sigma_2
\end{eqnarray*}
Note that $\M$ does not satisfying the standard ellipticity condition (see \cite[(3)]{tambue2016exponential}), so the PDE \eqref{conservation} is degenerated.

We will assume Dirichlet boundary condition in the entire domain.
\subsection{Finite volume method}

Let us consider the new domain $\Omega$ of study by truncating $D$ such that  $\Omega=I_x\times I_y\times [0,T]$ where $I_x=[0,x_{\max}]$ and $I_y=[0,y_{\max}]$.
In  the sequel  of  this work, the Black-Scholes partial differential equation \eqref{twoop} is considered  over the truncated domain $\Omega$. 
At $x=x_{\max}$ and $y=y_{\max}$, the linear boundary condition will be applied (see \cite{huang2006fitted}).
The intervals $I_x$ and $I_y$ will be subdivided into  $N$ part  in the following way (see \cite{huang2006fitted,huang2009convergence}) without loss the generality as irregular grids such as triangular grids can be used.  
\begin{eqnarray}
I_{x_i}=[x_{i-1};x_i], \,\,I_{y_j}=[y_{j-1};y_j]\quad i,j =1,...,N+1.
\end{eqnarray}
Let us set the mid-points $x_{i-\frac{1}{2}}$ and $y_{j-\frac{1}{2}}$ as follows
\begin{eqnarray}
x_{i-\frac{1}{2}}=\frac{x_{i-1}+x_i}{2}, \,\,\, y_{j-\frac{1}{2}}=\frac{y_{j-1}+y_j}{2} \qquad i,j =1,...,N,
\end{eqnarray}
with  $h_i=x_{i+\frac{1}{2}}-x_{i-\frac{1}{2}},~~~l_j=y_{j+\frac{1}{2}}-y_{j-\frac{1}{2}}$~~~~~and
\begin{eqnarray*}
x_{-\frac{1}{2}}=x_0=0, \qquad  x_{N+\frac{3}{2}}=x_{N+1}=x_{\max},\,  y_{-\frac{1}{2}}=y_0=0 \,, \;\,y_{N+\frac{3}{2}}=y_{N+1}=y_{\max}.
\end{eqnarray*}
For $i,j=1,\ldots,N$,  we denote by $\C_{ij}=[x_{i-\frac{1}{2}};x_{i+\frac{1}{2}}]\times[y_{j-\frac{1}{2}};y_{j+\frac{1}{2}}]$  a  control volume associated to our subdivision. 
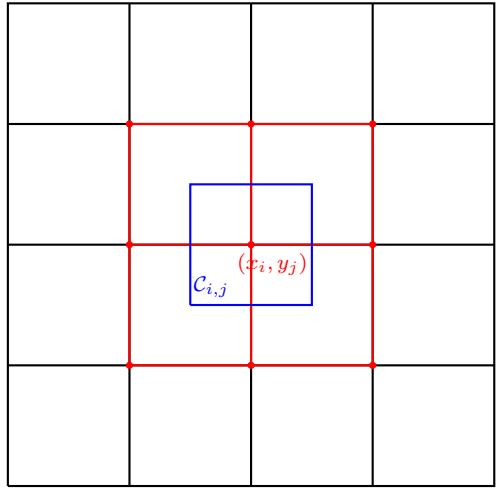
\begin{figure}[!h]
\centering
\begin{tikzpicture}[scale=0.4]
\draw[black,thick] (0,0)--(0,16)--(16,16)--(16,0)--(0,0);
\draw[black,thick] (4,0)--(4,16);
\draw[black,thick] (8,0)--(8,16);
\draw[black,thick] (12,0)--(12,16);
\draw[black,thick] (0,4)--(16,4);
\draw[black,thick] (0,8)--(16,8);
\draw[black,thick] (0,12)--(16,12);
\draw[red,thick] (4,4)--(4,12)--(12,12)--(12,4)--(4,4);
\draw[red,thick] (8,4)--(8,12);
\draw[red,thick] (4,8)--(12,8);
\draw[blue,thick] (6,6)--(6,10)--(10,10)--(10,6)--(6,6);
\draw[red,fill=red] (4,4) circle (0.1);
\draw[red,fill=red] (4,8) circle (0.1);
\draw[red,fill=red] (4,12) circle (0.1);
\draw[red,fill=red] (8,4) circle (0.1);
\draw[red,fill=red] (8,8) circle (0.1);
\draw[red,fill=red] (8,12) circle (0.1);
\draw[red,fill=red] (12,4) circle (0.1);
\draw[red,fill=red] (12,8) circle (0.1);
\draw[red,fill=red] (12,12) circle (0.1);
\node[below,red] at (8.7,8){$(x_i,y_j)$};
\node[above,blue] at (6.7,6){$\C_{i,j}$};
\end{tikzpicture}
\caption{The control volume $\C_{i,j}$}
\end{figure}
Note  that the control volume $\C_{ij}$ is the area surrounding the grid point $(x_i,y_j)$. 
Our goal  is to approximate  the option  function $U$ at $(x_i,y_j)$  \footnote{center of the control volume $\C_{i,j}$}  by a function denoted $\U$.
%
 The matrix $\M$ in \eqref{conservation} will be replaced by its average value  within each control  volume

\begin{equation}
\M^{ij}=\frac{1}{\mathrm {meas}(\C_{i,j})}\int_{\C_{i,j}}\M dxdy,\,\,\,  i,j=1,...,N.
\end{equation}
where $~\mathrm{meas}(\C_{ij})$ is the measure of $\C_{ij}$.
Thereby, we have  
\begin{equation*}
M^{i,j}= \left[\begin{array}{lcr}
\frac{\sigma_1^2}{6}\frac{x_{i+\frac{1}{2}}^3-x_{i-\frac{1}{2}}^3}{x_{i+\frac{1}{2}}-x_{i-\frac{1}{2}}}&    &\frac{\rho\sigma_1\sigma_2}{8}(x_{i+\frac{1}{2}}+x_{i-\frac{1}{2}})(y_{j+\frac{1}{2}}+y_{j-\frac{1}{2}}) \\
                   &        &  \\
\frac{\rho\sigma_1\sigma_2}{8}(x_{i+\frac{1}{2}}+x_{i-\frac{1}{2}})(y_{j+\frac{1}{2}}+y_{j-\frac{1}{2}})  &    &  \frac{\sigma_2^2}{6}\frac{y_{j+\frac{1}{2}}^3 -y_{j-\frac{1}{2}}^3}{y_{j+\frac{1}{2}}-y_{j-\frac{1}{2}}}
\end{array}\right].
\end{equation*}
\\
\\
Now let us consider the divergence form given in \eqref{conservation}. Following the finite volume method's principle,
we integrate the partial differential equation \eqref{conservation} over each control volume $\C_{ij}$ and we have
\begin{align}
\label{eqfinvol}
\int_{\C_{ij}}\frac{\partial{U}}{\partial{\tau}}d\C=\int_{\C_{ij}}\nabla \cdot(\M\nabla U) d\C+\int_{\C_{ij}}\nabla (f U)d\C
+\int_{\C_{ij}}\lambda U d\C.
\end{align}

The next section will be dedicated to spatial discretization of equation \eqref{eqfinvol}.
For the term in the  left hand side of  \eqref{eqfinvol} and for the last term in  its right hand side, we use the mid-point quadrature rule for their approximations.  More precisely 
\begin{eqnarray}
\int_{{\mathcal{C}}_{ij}}\frac{\partial{U}}{\partial{\tau}}d\C \approx  \mathrm {meas}(\C_{ij})\frac{d\U}{d\tau}(x_i,y_j,\tau)
\end{eqnarray}
\begin{eqnarray}
\label{linearterm}
\int_{{\mathcal{C}}_{ij}}\lambda U d\C \approx meas(\C_{ij})\lambda \U(x_i,y_j,\tau).
\end{eqnarray}
The diffusion term 
\begin{equation}
\label{diffusionterm}
\int_{\C_{ij}}\nabla \cdot(\M\nabla \U) d\C
\end{equation}
of \eqref{eqfinvol} will be approximated using the \textbf{Multi-point flux approximation} (MPFA) method or  our novel \textbf{fitted Multi-point flux approximation}.
More details will be given in the next section. Besides, the convection term 
\begin{equation}
\label{convectionterm}
\int_{\C_{ij}}\nabla (f \U)d\C
\end{equation}
of \eqref{eqfinvol} will be approximated using the upwind methods (first or second order).
Note that the standard two -point flux approximation in \cite{tambue2016exponential} can only be consistent in the approximation of  \eqref{diffusionterm} if and  only if the grid is $\M-$orthogonal. 
\section{Space discretization}
The spatial discretization of \eqref{conservation} consists of approximating all terms in \eqref{eqfinvol} over the control volumes of the study domain. 
\subsection{Discretization of the diffusion term}
Let us start by applying the divergence theorem to the diffusion term \eqref{diffusionterm} as follows, for $~~i,j=1,...,N$ 
\begin{equation}
\label{diffterm-disc}
\F^{ij}=\int_{\C_{ij}}\nabla \cdot(\M^{ij}\nabla \U)=\int_{\partial \C_{ij}}(\M^{ij}\nabla \U)\cdot\vec{n}d\partial\C
\end{equation}

where $\vec{n}$ is the outward vector from the control volume.\\
\\
Now, we can apply the so-called $\textbf{Multi-Point Flux Aprroximation(MPFA)}$  to approximate the integral defined in \eqref{diffterm-disc}.

\subsubsection{Multi-Point Flux Approximation (MPFA) method }  
There exists several types of Multi-Point Flux Approximation methods. The most known of MPFA methods are the O-method and the L-method. 
In our study, we focus on the O-method because it is the classical MPFA method and it is more intuitive comparing to the L-method which is fairly new and less intuitive (see \cite{aavatsmark2002introduction}).
Here, we follow  the description of the O-method developed  by \cite{aavatsmark2002introduction}.\\
We will start by  giving  an approximation of the gradient in the integral expression \eqref{diffterm-disc}.
\begin{itemize}
\item[]
Let us consider a triangle $x_1x_2x_3$, $\nu_i$ the outer normal vector of the edge located opposite of vertex $x_i,~i=1,2,3$ and $f$ a linear function over this triangle (see \figref{fig:triangle}).
The length of $\nu_i$ is equal to the length of the edge to which it is normal.
\begin{figure}[hbtp]
\centering
\subfloat[Triangle and corresponding  normal vectors ]{\begin{tikzpicture}[scale=0.8]
\draw[blue,thick] (2,2)--(3,5)--(6,1)--(2,2);
\draw[blue,thick,->] (2.5,3.5)--(0.5,12.5/3);
\draw[blue,thick,->] (4.5,3)--(7.5,5.6);
\node[above] at (2.8,3.2){$\bar{x}_A$};
\node[above] at (4.2,2.6){$\bar{x}_B$};
\node[above] at (1.8,1.5){$x_3$};
\node[above] at (6.2,0.5){$x_2$};
\node[above] at (3,5){$x_1$};
\node[above] at (0.5,4.1){$\nu_2$};
\node[above] at (7.5,5.6){$\nu_3$};
\end{tikzpicture}\label{fig:triangle}}
\qquad
\subfloat[A triangle in a control volume]{\begin{tikzpicture}[scale=0.6]
		\draw[blue,dotted, thick] (4,0)--(4,8);
		\draw[blue,dotted,thick] (0,4)--(8,4);
		\node[above,blue] at (2.3,1.3){$\mathcal{C}_{i,j}$};
		\draw[blue,thick] (2,2)--(6,2)--(6,6)--(2,6)--(2,2);
		\draw[thick] (0,0)--(0,8);
		\draw[thick] (0,0)--(8,0);
		\draw[thick] (8,0)--(8,8);
		\draw[thick] (0,8)--(8,8);
		\draw[->,thick] (4,5)--(6,5);
		\draw[->,thick] (5,4)--(5,6);
		\draw[red,thick] (4,4)--(4,6)--(6,4)--(4,4);
		\draw[blue,fill=blue] (6,6) circle (0.1);
		\draw[blue,fill=blue] (2,6) circle (0.1);
		\draw[blue,fill=blue] (6,2) circle (0.1);
		\draw[blue,fill=blue] (2,2) circle (0.1);
		\draw[red,fill=red] (4,6) circle (0.1);
		\draw[red,fill=red] (6,4) circle (0.1);
		\draw[red,fill=red] (4,4) circle (0.1);
		\node[below] at (4,4){$x_{1}$};
		\node[above] at (3.9,6.1){$\bar{x}_{2}$};
		\node[above] at (6.4,5){$\omega_1$};
		\node[below] at (5,6.8){$\omega_{2}$};
		\node[above] at (6.3,4){$\bar{x}_{1}$};
		\end{tikzpicture}\label{fig:triangle in CV}}
	\caption{}
\end{figure}
The gradient expression of the function $f$ in the triangle may be written in the form
\begin{equation}
\label{grad-tri}
\nabla f=-\frac{1}{2\mathcal{A}}\left[\Big(f(x_2)-f(x_1)\Big)\nu_2+\Big(f(x_3)-f(x_1)\Big)\nu_3\right]
\end{equation}

where $\mathcal{A}$ is the area of the triangle.

\newpage 
Thereby, assuming that our solution $\U$ is linear over  the   control volume $\C_{ij}$ with center $x_1(x_{i},y_{j})$, 
and applying \eqref{grad-tri} in the triangle $x_1\bar{x}_1\bar{x}_2$ (see \figref{fig:triangle in CV}), we have
  \begin{equation}
\label{grad-cv} 
\nabla \U=\frac{1}{2 \mathcal{A}}\left[(\bar{\U}_1-\U_{ij})\omega_1+(\bar{\U}_2-\U_{ij})\omega_2\right]
\end{equation}
 where $\U_{ij}=\U(x_1)=\U(x_{i},y_{j}),\bar{\U}_1=\U(\bar{x}_1),\bar{\U}_2=\U(\bar{x}_2)$ and the vectors $\omega_1$ and $\omega_2$ are respectively inner normal vector to the edge $x_1\bar{x}_1$ 
 and $x_1\bar{x}_2$ with the same length with those vectors, and  $\mathcal{A}$ is the area of the triangle $x_1\bar{x}_2\bar{x}_2$.\\
Let us called \textbf{interaction volume} $\R_{ij}$ a cell grid defined as follows
\begin{equation}
~~~\R_{ij}=[x_{i-1};x_i]\times[y_{j-1};y_j],\,\,i,j=1,\ldots,N+1.
\end{equation}

We may notice that an  interaction volume $\R_{ij}$ is covering an area  in the intersection of the control volumes $\C_{i-1,j-1},\C_{i-1,j},\C_{i,j-1}$ and $\C_{ij}$.
Here, we follow closely \cite{aavatsmark2007multipoint}.
 \item[]
 We denote respectively by $x_1(x_{i-1},y_{j-1}),x_2(x_{i},y_{j-1}),x_3(x_{i-1},y_{j})$
 ~and~$ x_4(x_{i},y_{j})$ the centre of the control volume $\C_{i-1,j-1},\C_{i,j-1},\C_{i-1,j}$ and $\C_{i,j}$. We denote also by $\bar{x}_1,\bar{x}_2,\bar{x}_3~and~\bar{x}_4$
 the midpoints of the segment $x_1x_2$, $x_3x_4,x_1x_3~ \text{and}~x_2x_4$
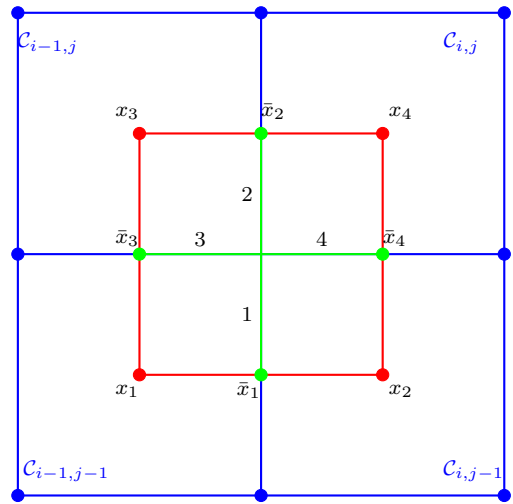
\begin{figure}[!h]
\begin{center}
\begin{tikzpicture}[scale=0.8]
\draw[blue,thick] (4,0)--(4,8);
\draw[blue,thick] (0,4)--(8,4);
\node[above,blue] at (0.8,0.1){$\mathcal{C}_{i-1,j-1}$};
\node[above,blue] at (7.5,0.1){$\mathcal{C}_{i,j-1}$};
\node[below,blue] at (7.3,7.8){$\mathcal{C}_{i,j}$};
\node[below,blue] at (0.5,7.8){$\mathcal{C}_{i-1,j}$};
\draw[red,thick] (2,2)--(6,2);
\draw[red,thick] (2,6)--(6,6);
\draw[red,thick] (2,2)--(2,6);
\draw[red,thick] (6,2)--(6,6);
\draw[blue,thick] (0,0)--(0,8);
\draw[blue,thick] (0,0)--(8,0);
\draw[blue,thick] (8,0)--(8,8);
\draw[blue,thick] (0,8)--(8,8);
\draw[green,thick] (2,4)--(6,4);
\draw[green,thick] (4,2)--(4,6) ;
\draw[red,fill=red] (6,6) circle (0.1);
\draw[red,fill=red] (2,6) circle (0.1);
\draw[red,fill=red] (6,2) circle (0.1);
\draw[red,fill=red] (2,2) circle (0.1);
\draw[blue,fill=blue] (0,0) circle (0.1);
\draw[blue,fill=blue] (0,4) circle (0.1);
\draw[blue,fill=blue] (0,8) circle (0.1);
\draw[blue,fill=blue] (4,0) circle (0.1);
\draw[blue,fill=blue] (4,8) circle (0.1);
\draw[blue,fill=blue] (8,0) circle (0.1);
\draw[blue,fill=blue] (8,4) circle (0.1);
\draw[blue,fill=blue] (8,8) circle (0.1);
\node[above] at (3,4){$3$};
\node[left] at (4,5){$2$};
\node[left] at (4,3){$1$};
\node[above] at (5,4){$4$};
\draw[green,fill=green] (4,6) circle (0.1);
\draw[green,fill=green] (6,4) circle (0.1);
\draw[green,fill=green] (4,2) circle (0.1);
\draw[green,fill=green] (2,4) circle (0.1);
\node[above] at (4.2,6.1){$\bar{x}_{2}$};
\node[above] at (1.8,4){$\bar{x}_{3}$};
\node[below] at (3.8,2){$\bar{x}_{1}$};
\node[above] at (6.2,4){$\bar{x}_{4}$};
\node[above] at (1.8,1.5){$x_1$};
\node[above] at (6.3,1.5){$x_2$};
\node[above] at (1.8,6.1){$x_3$};
\node[above] at (6.3,6.1){$x_4$};
\end{tikzpicture}
\end{center}
\caption{Interaction volume}
\label{fig:intervol}
\end{figure}

Our goal in an interaction volume is to compute the flux through the half edges $1,2,3~and~4$ inside the interaction volume (see \figref{fig:intervol}).
The flux through the half edge p seen from the centre $x_1=(x_{i-1},y_{j-1})$ of the control volume $\C_{i-1,j-1}$ is denoted $f_p^{i-1,j-1}$. By using the expression \eqref{diffterm-disc}, we have
\begin{equation}
\label{flux-expr}
f_p^{i-1,j-1}=\Gamma_p\vec{n}_p^T\M^{i-1,j-1}\nabla \U
\end{equation}
where $\Gamma_p$ is the length of half edge p, $\vec{n}_p$ is the outward unit normal vector to the half edge p. It is convenient to let $\vec{n}_p$ point in the direction of increasing 
global cell indices. In that case, we have two kinds of  inner normal vectors. The vertical ones denoted $\omega_1$ and the horizontal ones denoted $\omega_2$.\\
By considering the triangle $x_1\bar{x}_1\bar{x}_3$ (\figref{fig:intervol}) in the control volume $\C_{i-1,j-1}$, using the expression of gradient \eqref{grad-cv} and the flux expression 
\eqref{flux-expr}, we have for $i,j=1,...,N$
\begin{equation}
\label{flux-vec}
\left[\begin{array}{c}
f_1^{i-1,j-1} \\   \\ f_3^{i-1,j-1} 
\end{array}\right] = G^{i-1,j-1}\left[\begin{array}{c}
\bar{\U}_1-\U_{i-1,j-1}  \\ \\ \bar{\U}_3-\U_{i-1,j-1}
\end{array} \right]
\end{equation}  
with
\begin{equation*}
G^{i-1,j-1}=\begin{bmatrix}
\Gamma_1n_1^T M^{i-1,j-1} \omega_1  &  &   &  & \Gamma_1n_1^T M^{i-1,j-1} \omega_2 \\
    & &   &   & \\
 \Gamma_2n_1^T M^{i-1,j-1} \omega_2 &   &  &  & \Gamma_2n_2^T M^{i-1,j-1} \omega_2
\end{bmatrix}
\end{equation*}
%

By applying \eqref{flux-vec} in the triangles $x_2\bar{x}_1\bar{x}_4,x_3\bar{x}_2\bar{x}_3$ and $x_4\bar{x}_4\bar{x}_2$ (see  \figref{fig:intervol}), we have
\begin{eqnarray}
\label{flux-vec1}
\left[\begin{array}{c}
f_1^{i,j-1} \\   \\ f_4^{i,j-1} 
\end{array}\right] = G^{i,j-1}\left[\begin{array}{c}
\U_{i,j-1}-\bar{\U}_1 \\ \\\bar{\U}_4-\U_{i,j-1}
\end{array} \right]~~~~~~~~~~~~~\left[\begin{array}{c}
f_2^{i-1,j} \\   \\ f_3^{i-1,j} 
\end{array}\right] = G^{i-1,j}\left[\begin{array}{c}
\bar{\U}_2-\U_{i-1,j}  \\ \\ \U_{i-1,j}-\bar{\U}_3
\end{array} \right] \nonumber\\
\nonumber \\
\\
\nonumber \\
\left[\begin{array}{c}
f_2^{ij} \\   \\ f_4^{ij} 
\end{array}\right] = G^{ij}\left[\begin{array}{c}
\U_{ij}-\bar{\U}_2  \\ \\ \U_{ij}-\bar{\U}_4
\end{array} \right] \nonumber ~~~~~~~~~~~~~~~~~~~~~~~~~~~~~~~~~
\end{eqnarray}
Since the flux through an edge is continuous, from \eqref{flux-vec} and \eqref{flux-vec1} we have 
\begin{equation}
\begin{array}{ccccc}
f_1 & = & f_1^{i-1,j-1} &  = &  f_1^{i-1,j} \\
    &   &  &     &  \\
 f_2 &  = &  f_2^{ij} & = & f_2^{i-1,j} \\
   &  &    &  &    \\
 f_3 &   =   &  f_3^{i-1,j} & = & f_3^{i-1,j-1} \\
  &   &    &    &    \\
 f_4  &   =   &      f_4^{i,j-1}  & =  & f_4^{ij}.  
\end{array}
\end{equation}
It follows that 
\begin{eqnarray}
\label{flux-cont}
f_1 & = & g_{11}^{i-1,j-1}(\bar{\U}_1-\U_{i-1,j-1})+g_{12}^{i-1,j-1}(\bar{\U}_3-\U_{i-1,j-1})
  =  -g_{11}^{i,j-1}(\bar{\U_1}-\U_{i,j-1})+g_{12}^{i,j-1}(\bar{\U}_4-\U_{i,j-1}) \nonumber \\
&    &      \nonumber  \\
f_2 & = & -g_{11}^{ij}
(\bar{\U}_2-\U_{ij})-g_{12}^{ij}(\bar{\U}_4-\U_{ij})
 =  g_{11}^{i-1,j}(\bar{\U}_2-\U_{i-1,j})-g_{12}^{i-1,j}(\bar{\U}_3-\U_{i-1,j}) \nonumber \\
 &        &    \\
 f_3 & = & g_{21}^{i-1,j}(\bar{\U}_2-\U_{i-1,j})-g_{22}^{i-1,j}(\bar{\U}_3-\U_{i-1,j})
   =  g_{21}^{i-1,j-1}(\bar{\U}_1-\U_{i-1,j-1})+g_{22}^{i-1,j-1}(\bar{\U}_3-\U_{i-1,j-1})\nonumber \\
   &         &  \nonumber  \\
f_4 & = & -g_{21}^{i,j-1}(\bar{\U}_1-\U_{i,j-1})+g_{22}^{i,j-1}(\bar{\U}_4-\U_{i,j-1})
  =  -g_{21}^{ij}(\bar{\U}_2-\U_{ij})-g_{22}^{ij}(\bar{\U}_4-\U_{ij}) \nonumber
\end{eqnarray}

Let us set 
\begin{equation}
f=\left[\begin{array}{c} f_1 \\ f_2 \\  f_3\\ f_4\end{array}\right],~~~~~~~~~~\U=\left[\begin{array}{c} \U_{i-1,j-1} \\   \U_{i,j-1} \\  \U_{i-1,j} \\ 
\U_{ij}\end{array}\right],~~~~~~~\V=\left[\begin{array}{c} \bar{\U}_1 \\  \bar{\U}_2 \\  \bar{\U}_3 \\  \bar{\U}_4\end{array}\right]
\end{equation}
The equation \eqref{flux-cont} allows to have 
\begin{equation}
\label{flux-eq1}
f=C^{ij}\V+F^{ij}\U
\end{equation}
where
\begin{eqnarray*}
 C^{ij} & = & \left[\begin{array}{ccccccc}
g_{11}^{i-1,j-1} &  & 0 &  & g_{12}^{i-1,j-1} & & 0 \\
&   & &  & &  & \\
0  &  & -g_{11}^{ij}  &  & 0 & & -g_{12}^{ij} \\
&  &  &  &   & & \\
0  &  & g_{21}^{i-1,j} &  & -g_{22}^{i-1,j} &  &  0 \\
&  &  &  &   &  & \\
-g_{21}^{i,j-1}  &  & 0 &  & 0 &  &   g_{22}^{i,j-1}
\end{array}\right]
\end{eqnarray*}

\begin{eqnarray*}
F^{ij}=\left[\begin{array}{ccccccc}
-g_{11}^{i-1,j-1}-g_{12}^{i-1,j-1} &  & 0 & &  0 & &  0 \\
&  &   &   &   &   &    \\
0 & &  0  &   &  0 &  &  g_{11}^{ij}+g_{12}^{ij}\\
&   &   &   &   &  & \\
0 &   &  0 &  &  -g_{21}^{i-1,j}+g_{22}^{i-1,j} &   &  0 \\\\
&   &   &   &   &   & \\
0 &   &  g_{21}^{i,j-1}-g_{22}^{i,j-1} &   &  0   &  &    0 \end{array}\right]
\end{eqnarray*}
From \eqref{flux-cont}, we can also have 
\begin{eqnarray}
\label{flux-eq2}
A^{ij}\V=B^{ij}\U
\end{eqnarray}
where 
\begin{eqnarray*}
A^{ij} & = &\left[\begin{array}{ccccccc}
g_{11}^{i-1,j-1}+g_{11}^{i,j-1}  & & 0 & &  g_{12}^{i-1,j-1} & & -g_{12}^{i,j-1}\\
&   &   &   &    &    & \\
0 &   &  -g_{11}^{ij}-g_{11}^{i-1,j} &  & g_{12}^{i-1,j} &   & -g_{12}^{ij}\\
&    &    &   &   &    & \\
-g_{21}^{i-1,j-1} &    & g_{21}^{i-1,j} &   & -g_{22}^{i-1,j}-g_{22}^{i-1,j-1} &   & 0 \\
&    &    &    &    &    & \\
-g_{21}^{i,j-1} &    & g_{21}^{ij} &   & 0 &    & g_{22}^{i,j-1}+g_{22}^{ij}
\end{array}\right]
\end{eqnarray*}
\begin{eqnarray*}
B^{ij} & = & \left[\begin{array}{ccccccc}
g_{11}^{i-1,j-1}+g_{12}^{i-1,j-1} &   &  g_{11}^{i,j-1}-g_{12}^{i,j-1} &  &  0 &  &  0\\
&    &    &    &    &    & \\
0 &   &  0 &   &  -g_{11}^{i-1,j}+g_{12}^{i-1,j}  &    & -g_{11}^{ij}-g_{12}^{ij}\\
&     &    &    &    &    & \\
-g_{21}^{i-1,j-1}-g_{22}^{i-1,j-1} & &  0 & &  g_{21}^{i-1,j}-g_{22}^{i-1,j} &    & 0 \\
&    &    &   &    &    & \\
0 &   &  -g_{21}^{i,j-1}+g_{22}^{i,j-1} &   &  0 &   &  g_{21}^{ij}+g_{22}^{ij}\end{array}\right]
\end{eqnarray*}
Thereby, $\V$ can be eliminated from \eqref{flux-eq1} by solving \eqref{flux-eq2} with respect to $\V$. This gives the following  
the expression of the flux through the 4 half edges inside the interaction  volume  $\R_{ij}$ 
\begin{eqnarray}
\label{flux-trans}
f=T^{ij}\U, \;\;\;\,\,\,\, i,j=1,...,N+1.
\end{eqnarray}
where 
\begin{equation}
\label{trans}
T^{ij}=C^{ij}\left[A^{ij}\right]^{-1}B^{ij}+F^{ij}
\end{equation}
$T^{ij}$ is called transmissibility matrix of the interaction volume $\R_{ij}$.\\
 From \eqref{flux-trans}, we are now able to get the flux through the half edges 1,2,3 and 4 inside the interaction volume $\R_{ij}$.\\
Let us recall that  to approximate  the integral in \eqref{diffterm-disc}, we need to compute the flux through the edges on a control volume $\C_{ij}$. 
We might notice that we need four interaction volume with centres the four vertices of the control volumes in order to cover all the edges of the considered control volume (see \figref{fig:intervol1}).
\begin{figure}[!h]
\centering
\begin{tikzpicture}[scale=0.5]
\draw[black,thick] (0,0)--(0,16)--(16,16)--(16,0)--(0,0);
\draw[black,thick] (4,0)--(4,16);
\draw[black,thick] (8,0)--(8,16);
\draw[black,thick] (12,0)--(12,16);
\draw[black,thick] (0,4)--(16,4);
\draw[black,thick] (0,8)--(16,8);
\draw[black,thick] (0,12)--(16,12);
\draw[red,thick] (4,4)--(4,12)--(12,12)--(12,4)--(4,4);
\draw[red,thick] (8,4)--(8,12);
\draw[red,thick] (4,8)--(12,8);
\draw[blue,dotted] (10,6)--(6,6)--(6,10)--(10,10);
\draw[black,thick]  (10,6)--(10,10);
\node[below,black] at (10.2,8){$\E$};
\draw[red,fill=red] (4,4) circle (0.1);
\draw[red,fill=red] (4,8) circle (0.1);
\draw[red,fill=red] (4,12) circle (0.1);
\draw[red,fill=red] (8,4) circle (0.1);
\draw[red,fill=red] (8,8) circle (0.1);
\draw[red,fill=red] (8,12) circle (0.1);
\draw[red,fill=red] (12,4) circle (0.1);
\draw[red,fill=red] (12,8) circle (0.1);
\draw[red,fill=red] (12,12) circle (0.1);
\node[below,red] at (8.7,8){$(x_i,y_j)$};
\node[above,blue] at (6.7,6){$\C_{ij}$};
\node[above] at (9.8,6.2){$2$};
\node[above] at (9.8,8.2){$1$};
\draw[black,fill=black] (10,8) circle (0.1);
\node[above,red] at (4.5,4.1){$\R_{ij}$};
\node[above,red] at (13.5,4.3){$\R_{i+1,j}$};
\node[above,red] at (13.7,8.7){$\R_{i+1,j+1}$};
\node[above,red] at (4.5,8.7){$\R_{i,j+1}$};
\end{tikzpicture}
\caption{}
\label{fig:intervol1}
\end{figure}

\newpage
For the volume control $\C_{ij}$, we denote by ${}_{\E}f_{l}^{ij}$ the flux through lower half eastern edge, by ${}_{\E}f_{u}^{ij}$ the flux through the upper 
half eastern edge. The flux  ${}_{\E}f^{ij}$  through the east edge of the control volume $\C_{ij}$ is calculated as follows:
The lower half eastern edge is contained in the interaction volume $\R_{i+1,j}$  and it is in position 2 in the interaction of volume (see \figref{fig:intervol1}).
So by using \eqref{flux-trans} we have:

\begin{equation*}
{}_{\E}f_{l}^{ij}=T_{21}^{i+1,j}\U_{i,j-1}+T_{22}^{i+1,j}\U_{i+1,j-1}
+T_{23}^{i+1,j}\U_{ij}+T_{24}^{i+1,j}\U_{i+1,j}.
\end{equation*}
Similarly, the upper half eastern edge  is contained in the interaction volume $\R_{i+1,j+1}$  and it is in position 1 in the interaction volume. So by using \eqref{flux-trans} we have: 
\begin{equation*}
{}_{\E}f_u^{ij}=T_{11}^{i+1,j+1}\U_{ij}+T_{12}^{i+1,j+1}\U_{i+1,j}+T_{13}^{i+1,j+1}\U_{i,j+1}
+T_{14}^{i+1,j+1}\U_{i+1,j+1}.
\end{equation*}

Finally the flux through the east edge of the control volume $\C_{i+1,j+1}$ will be the addition of ${}_{\E}f_{l}^{ij}$ and  ${}_{\E}f_u^{ij}$. Thereby we have
\begin{eqnarray*}
{}_{\E}f_{}^{ij} & = & {}_{\E}f_{l}^{ij}+{}_{\E}f_u^{ij}  \\
    &   & \\
    & = &    T_{21}^{i+1,j}\U_{i,j-1}+T_{22}^{i+1,j}\U_{i+1,j-1}
+T_{23}^{i+1,j}\U_{ij}+T_{24}^{i+1,j}\U_{i+1,j}+
 T_{11}^{i+1,j+1}\U_{ij}\\
 &   & \\
&  & +T_{12}^{i+1,j+1}\U_{i+1,j}+T_{13}^{i,j}\U_{i,j+1}+T_{14}^{i,j}\U_{i+1,j+1} \\
&   & \\
{}_{\E}f_{}^{ij}  & = & (T_{11}^{i+1,j+1}+T_{23}^{i+1,j})\U_{ij}+(T_{12}^{i+1,j+1}+T_{24}^{i+1,j})\U_{i+1,j}+
T_{14}^{i+1,j+1}\U_{i+1,j+1}\\
& &\\
&   &  +T_{13}^{i+1,j+1}\U_{i,j+1}+T_{21}^{i+1,j}\U_{i,j-1}+T_{22}^{i+1,j}\U_{i+1,j-1}.
\end{eqnarray*}
Similarly, we compute the flux through the northern, western and southern edge of the control volume $\C_{ij}$. 
Afterwards, we sum up the flux through the 4 edges of the control to get the outflux $\F^{ij}$ through the edges of the control volume $\mathcal{C}_{ij}$.
Therefore we have for $i,j=1,...,N$ 
\begin{eqnarray}
\label{flux-mpfa}
\F^{ij}& = & a_{ij}\U_{ij}+b_{ij}\U_{i+1,j}+c_{ij}\U_{i+1,j+1}+d_{ij}\U_{i,j+1}
+e_{ij}\U_{i-1,j+1}+\alpha_{ij}\U_{i-1,j}+\beta_{ij}\U_{i-1,j-1} 
 \nonumber \\
&  & +\gamma_{ij}\U_{i,j-1}+\lambda_{ij}\U_{i+1,j-1}.
\end{eqnarray} 
where 
\begin{eqnarray*}
	& &a_{ij}=T_{11}^{i+1,j+1}+T_{23}^{i+1,j}+T_{31}^{i+1,j+1}+T_{42}^{i,j+1}-T_{12}^{i,j+1}-T_{24}^{ij}-T_{33}^{i+1,j}-T_{44}^{ij};\\
	& & \\
	& & b_{ij}=T_{12}^{i+1,j+1}+T_{24}^{i+1,j}+T_{32}^{i+1,j+1}-T_{34}^{i+1,j}\\
	& & \\
	& &c_{ij}=T_{14}^{i+1,j+1}+T_{34}^{i+1,j+1}; d_{ij}=T_{13}^{i+1,j+1}+T_{33}^{i+1,j+1}+T_{44}^{i,j+1}-T_{14}^{i,j+1};
	e_{ij}=T_{43}^{i,j+1}-T_{13}^{i,j+1};\\
	& & \\
	& &  \alpha_{ij}=T_{41}^{i,j+1}-T_{11}^{i,j+1}-T_{23}^{ij}-T_{43}^{ij}; 
	\beta_{ij}=-T_{21}^{ij}-T_{41}^{ij};\\
	& & \\
	& & \gamma_{ij}=T_{21}^{i+1,j}-T_{22}^{ij}-T_{31}^{i+1,j}-T_{42}^{ij};\\
	& & \\
	& & \lambda_{ij}=T_{22}^{i+1,j}-T_{32}^{i+1,j}.
\end{eqnarray*}

Let us  notice that for the control volumes   near to the boundary of the our domain,  some terms from the boundary conditions will be involved in \eqref{flux-mpfa} .\\
 Hence \eqref{diffterm-disc} becomes 
\begin{equation}
\label{flux-mpfa-mat}
\F=A_{mp}\U+F_{mp}
\end{equation} 
where $A_{mp}$ is a $N^2\times N^2$ matrix and 
\begin{equation*}
\F=\begin{bmatrix}
\F_{11}\\
\F_{12}\\
\vdots\\
\F_{1N}\\
\F_{21}\\
\F_{22}\\
\vdots\\
\vdots\\
\F_{NN}
\end{bmatrix},~~~
\U=\begin{bmatrix}
\U_{11}\\
\U_{12}\\
\vdots\\
\U_{1N}\\
\U_{21}\\
\U_{22}\\
\vdots\\
\vdots\\
\U_{NN}
\end{bmatrix},~~~
A_{mp}=\begin{bmatrix}
W_1 & X_1 & 0_N & \ldots & \ldots & \ldots & \ldots & 0_N\\ 
Y_2 & W_2 & X_2 & \ddots &  &  &  & \vdots  \\ 
0_N & Y_3 & W_3 & X_3 & \ddots &  & & \vdots  \\ 
\vdots & \ddots & Y_4 & W_4 & X_4 & \ddots & & \vdots\\ 
\vdots &  & \ddots & \ddots & \ddots & \ddots & \ddots & \vdots \\ 
\vdots &  &  & \ddots& \ddots & \ddots & \ddots & 0_N \\ 
\vdots &  &  &  & \ddots & Y_{N-1} & W_{N-1} & X_{N-1} \\ 
0_N & \ldots & \ldots & \ldots & \ldots & 0_N & Y_N & W_N
\end{bmatrix}   
\end{equation*}
with
$0_{N}$ is $N\times N$ null matrix , $W_i,Y_i,X_i$ are tridiagonal matrices, and $F_{mp}$ is a $N^2$ vector coming from the boundary conditions.  
The  structure of the diffusion matrix $A_{mp}$ can be viewed in \figref{diffusion1}
\begin{figure}[hbtp]
\centering
\includegraphics[scale=0.5]{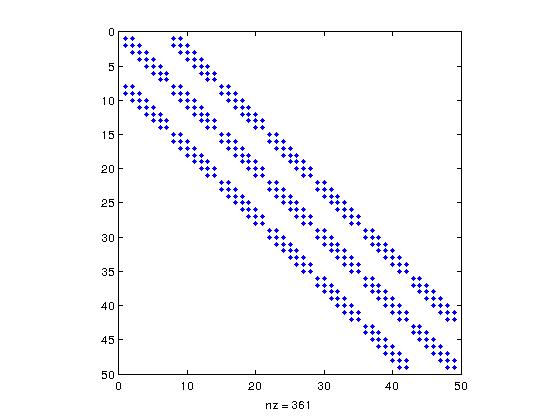}
\caption{Structure of diffusion matrix coming from standard MPFA}
\label{diffusion1}
\end{figure}
\end{itemize}
\newpage
\subsection{Discretization of the convection term}
In this section, the convection term 
\begin{equation*}
\int_{\C_{ij}}\nabla (f \U)d\C
\end{equation*}

with 
\begin{equation*}
f=\left(\begin{array}{c}
(r-\sigma_1^2-\frac{1}{2}\rho\sigma_1\sigma_2)x \\ \\ (r-\sigma_2^2-\frac{1}{2}\rho\sigma_1\sigma_2)y
\end{array}\right)=\left(\begin{array}{c}
p \\ \\ q
\end{array}\right)
\end{equation*}
will be approximated by  the upwind methods (first and second order).
\subsubsection{First order upwind}
The \textbf{first order upwind method} discussed by \cite[chapter 4.8]{leveque2004finite} or \cite{tambue2016exponential} will be applied to approximate the second term of \eqref{eqfinvol}.
Using the divergence theorem, we have  for $i,j=2,...,N$
\begin{equation}
I^{ij}=\int_{\C_{ij}}\nabla (f \U)d\C=\int_{\partial \C_{ij}}(f\cdot\U)\cdot\vec{n}d\partial \C.
\end{equation}
Note that $I^{ij}$ is calculated by summing up the flux through the edges of the control volume $\C_{ij}$.
The flux through an edge using the first order upwind will depend on the sign of $f\cdot\vec{n}$ on this edge. If the sign of $f\cdot\vec{n}$
is positive, $\U_{ij}$ will be used to approximate $\U$  in the expression $(f\cdot \vec{n} \U)$ otherwise we will use the value of $\U$ in other side of the edge. 
Note that an edge may be the interface of two control volumes.
By doing so, we  have  for $i,j=1,...,N$
\begin{eqnarray}
\label{flux-up1}
I^{ij} & =  & \epsilon_{ij}\U_{i-1,j}+\mu_{ij}\U_{i,j-1}+\Omega_{ij}\U_{ij}+
  \phi_{ij}\U_{i,j+1}+\Psi_{ij}\U_{i+1,j}.
\end{eqnarray}
where 
\begin{eqnarray*}
&   &\epsilon_{ij}=-l_jf_x^{i-1}\max(f_x^{i-1},0);~~~~~\mu_{ij}=-h_if_y^{j-1}\max(f_y^{j-1},0)\\
&   & \\
&   & \Omega_{ij}=l_j\Bigg(f_x^{i}\max(f_x^{i},0)- f_x^{i-1}\min(f_x^{i-1},0)\Bigg)
  +h_i\Bigg(f_y^j\max(f_y^j,0)-f_y^{j-1}\min(f_y^{j-1},0)\Bigg)\\
&  & \\
&   & \\
&   & \phi_{ij}=h_if_y^j\min(f_y^j,0);~~~~~~~~
\Psi_{ij}=l_jf_x^{i}\min(f_x^{i},0).
\end{eqnarray*}

with

\begin{eqnarray*}
f_x^i=(r-\sigma_1-\frac{1}{2}\rho\sigma_1\sigma_2)x_{i+1}~~~~~~~~~~~~~~~~f_y^j=(r-\sigma_2-\frac{1}{2}\rho\sigma_1\sigma_2)x_{j+1}
\end{eqnarray*}

Let us  notice that for the control volumes   near to the boundary of the our domain,  some terms from the boundary conditions will be involved in \eqref{flux-up1}.
Hence, \eqref{flux-up1} gives 
 \begin{equation}
 \label{flux-up1-mat}
I=A_{up} \U+F_{up}
\end{equation}
where  $A_{up}$ is a $N^2\times N^2$ matrix 
\begin{equation*}
I=\begin{bmatrix}
I^{11}\\
I^{12}\\
\vdots\\
I^{1N}\\
I^{21}\\
I^{22}\\
\vdots\\
\vdots\\
I^{NN}\\
\end{bmatrix},
~~~~~~\U=\begin{bmatrix}
\U_{11}\\
\U_{12}\\
\vdots\\
\U_{1N}\\
\U_{21}\\
\U_{22}\\
\vdots\\
\vdots\\
\U_{NN}\\
\end{bmatrix}, A_{up}=\begin{bmatrix}
H_1 & P_1 & 0_N & \ldots & \ldots & \ldots & 0_N \\ 
Q_2 & H_2 & P_2& \ddots & &   & \vdots \\
0_N & Q_3 & H_3 & P_3 & \ddots &   &   \vdots \\
\vdots & \ddots & \ddots & \ddots & \ddots & \ddots & \vdots \\
\vdots  &   & \ddots & Q_{N-2}  & H_{N-2} & P_{N-2} & 0_N \\
\vdots &  &  & \ddots & Q_{N-1} & H_{N-1} & P_{N-1}\\
 0_N & \dots   &   \ldots  & \ldots  &  0_N &  Q_N  & H_N
\end{bmatrix}
\end{equation*}
with $0_{N}$ is $N\times N$ null matrix, $H_i$ is a tridiagonal matrix, $P_i,Q_i$ are diagonal matrices and $F_{up}$ is a vector coming from the boundary conditions. 
%
Therefore, combining the MPFA method \eqref{flux-mpfa-mat} and the first order upwind \eqref{flux-up1-mat}, we have 
\begin{equation}
\label{mpfa-up1}
\frac{d\U}{d\tau}=A\U+F
\end{equation} 
with 
\begin{eqnarray*}
A=L^{-1}\Bigg(A_{mp}+A_{up}+A_L\Bigg)~~~~~~F=L^{-1}\Bigg(F_{mp}+F_{up}\Bigg)
\end{eqnarray*}

where $A_L$ is a diagonal matrix of size $N^2\times N^2$ coming from the discretisation of \eqref{linearterm}.
The diagonal elements of $A_L$ are $A_{ii}=h_il_i\lambda$ for $i=1,...,N$ with $\lambda$ given in \eqref{conservation}.
The matrix $L$ is also a diagonal matrix of size $N^2\times N^2$ whose diagonal elements are $L_{ii}=h_il_i$ for $i=1,\ldots,N$
\subsubsection{Upwind second order}
We start by applying the mid-quadrature rule as follows.
\begin{eqnarray}
\label{quad-upw2}
J^{ij}=\int_{\C_{ij}}\nabla (f\U) d\C& = &  mes(\C_{ij})\nabla (f\U)|_{(x_i,y_j)}\nonumber\\
                           &   &   \nonumber \\
                           &  = & (x_{i+\frac{1}{2}}-x_{i-\frac{1}{2}})(y_{j+\frac{1}{2}}-y_{j-\frac{1}{2}})\Bigg[p_i\frac{\partial \U_{ij}}{\partial x}+q_j\frac{\partial \U_{ij}}{\partial y}+\Bigg(\frac{\partial p_i}{\partial x}+\frac{\partial q_j}{\partial y}\Bigg)\U_{ij}\Bigg]\nonumber\\
                           &  & \nonumber \\
                           &  = &  h_il_j\Bigg[\Bigg(p_i\frac{\partial \U_{ij}}{\partial x}+q_j\frac{\partial \U_{ij}}{\partial y}\Bigg)+\omega\U_{ij}\Bigg],\;\;\,i,j=1,\ldots,N.
\end{eqnarray}

where $~~~~~p_i=(r-\sigma_1^2-\frac{1}{2}\rho\sigma_1\sigma_2)x_i,~~~q_j=(r-\sigma_2^2-\frac{1}{2}\rho\sigma_1\sigma_2)y_j$ 
and $\omega=2r-\sigma_1^2-\sigma_2^2-\rho\sigma_1\sigma_2$.
Let us use the second order upwind to approximate the first derivatives in \eqref{quad-upw2} at the point $(x_i,y_j)$.\\

\textbf{Approximation of the first derivative using a 3 points stencil}

Here, we want to express the first derivative $\frac{\partial \U_{ij}}{\partial x}$ in terms of $\U_{i+2,j},\U_{i+1,j}$ and $\U_{ij}$.  Set $h=\underset{1\leq i\leq N}{\max} h_i$.  Let us find $a,b$ and $c$ such that
\begin{equation}
\label{first-dev}
\frac{\partial \U_{ij}}{\partial x}=a\U_{i+2,j}+b\U_{i+1,j}+c\U_{ij}
\end{equation}

\begin{figure}[hbtp]
	\centering
    \includegraphics[scale=0.5]{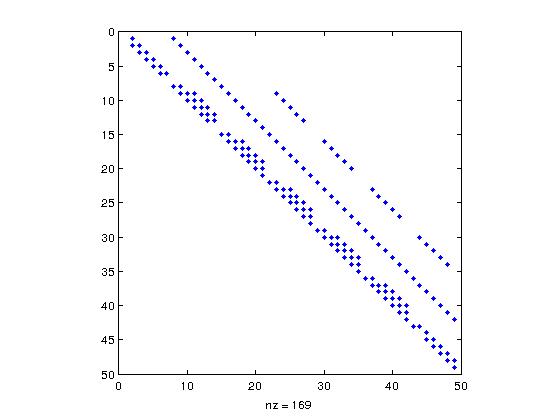}
    \caption{A structure of the  advection matrix using 2nd order upwind method.}
    \label{upw22}
\end{figure}

Thereby, using  a $2^{nd}$ order Taylor expansion at the point $(x_i,y_j)$ on $\U_{i+2,j}$ and $\U_{i+1,j}$, we have
\begin{eqnarray*}
\frac{\partial \U_{ij}}{\partial x}  &  =  & a\U_{i+2,j}+b\U_{i+1,j}+c\U_{ij}\\             
&     & \\
& = &  a\Bigg(\U_{ij}+(h_{i+1}+h_{i+2})\frac{\partial \U_{ij}}{\partial x}+ \frac{1}{2}(h_{i+1}+h_{i+2})^2\frac{\partial^2 \U_{ij}}{\partial x^2}+\mathcal{O}(h^3)\Bigg)+b\Bigg(\U_{ij}+h_{i+1}\frac{\partial \U_{ij}}{\partial x}+\frac{1}{2}h_{i+1}^2\frac{\partial^2 \U_{ij}}{\partial x^2}+\mathcal{O}(h^3)\Bigg)\\
&    & \\
&    & +c\U_{ij}.\\
&    & \\
\frac{\partial \U_{ij} }{\partial x}  &  =  & \Big(a+b+c\Big)\U_{ij}+\Bigg(a(h_{i+1}+h_{i+2})+bh_{i+1}\Bigg)\frac{\partial \U_{ij}}{\partial x}+\Bigg(\frac{1}{2}a\Big(h_{i+1}+h_{i+2}\Big)^2+\frac{1}{2}bh_{i+1}^2\Bigg)\frac{\partial^2 \U_{ij}}{\partial x^2}+\mathcal{O}(h^3).
\end{eqnarray*} 
By matching, we have
\begin{align}
\label{sys-first-dev}
\left\lbrace \begin{array}{l}
a+b+c=0\\
\\
a(h_{i+1}+h_{i+2})+bh_{i+1}~=1\\
\\
\frac{1}{2}a\Big(h_{i+1}+h_{i+2}\Big)^2+\frac{1}{2}bh_{i+1}^2=0
\end{array}\right.
\end{align}
Solving \eqref{sys-first-dev}, we have
\begin{equation}
a=-\frac{h_{i+1}}{h_{i+2}(h_{i+1}+h_{i+2})}~~~~~~~~~~~~~~~~b=\frac{h_{i+1}+h_{i+2}}{h_{i+1}h_{i+2}}~~~~~~~~~~~~~~~
c=\frac{h_{i+1}^2-\Big(h_{i+1}+h_{i+2}\Big)^2}{h_{i+1}+h_{i+2}}.
\end{equation}
Therefore we have
\begin{equation}
\label{first-dev-2nd}
\frac{\partial \U_{ij} }{\partial x} \approx \frac{-h_{i+1}^2\U_{i+2,j}+(h_{i+1}+h_{i+2})^2\U_{i+1,j}+(h_{i+1}^2-(h_{i+1}+h_{i+2})^2)\U_{ij}}{h_{i+1}h_{i+2}(h_{i+1}+h_{i+2})}.\\
\end{equation}\\
\textbf{Application to the $2^{nd}$ order upwind method on non uniform grids}\\
By analogy with the procedure to get the expression in \eqref{first-dev-2nd}, the term $p_i\frac{\partial\U_{ij}}{\partial x}$ is approximated as follows:
\begin{itemize}
\item[(i)] $p_i>0$ then 
\begin{equation*}
p_i\frac{\partial \U_{ij}}{\partial x} \approx p_i \frac{(h_{i+1}+h_{i+2})^2\U_{i+1,j}+\Big[h_{i+1}^2-(h_{i+1}+h_{i+2})^2\Big]\U_{ij}-h_{i+1}^2\U_{i+2,j}}{h_{i+1}h_{i+2}(h_{i+1}+h_{i+2})}
\end{equation*}
\item[(ii)] $p_i<0$ then
\begin{equation*}
p_i\frac{\partial \U_{ij}}{\partial x} \approx p_i \frac{-(h_i+h_{i-1})^2\U_{i-1,j}+\Big[(h_{i}+h_{i-1})^2-h_i^2\Big]\U_{ij}+h_i^2\U_{i-2,j}}{h_ih_{i-1}(h_i+h_{i-1})}
\end{equation*}
\end{itemize}
Similarly for the first derivative $\frac{\partial \U_{ij}}{\partial y}$, we have
\begin{itemize}
\item[(iii)] when $q_j>0$ then
\begin{equation*}
q_j\frac{\partial \U_{ij}}{\partial y} \approx q_j \frac{(l_{j+1}+l_{j+2})^2\U_{i,j+1}+\Big[l_{j+1}^2-(l_{j+1}+l_{j+2})^2\Big]\U_{ij}-l_{j+1}^2\U_{i,j+2}}{l_{j+1}l_{j+2}(l_{j+1}+l_{j+2})}
\end{equation*}
\item[(iv)] when $q_j<0$ 
\begin{equation*}
q_j\frac{\partial \U_{ij}}{\partial y} \approx q_j \frac{-(l_j+l_{j-1})^2\U_{i,j-1}+\Big[(l_j+l_{j-1})^2-l_j^2\Big]\U_{ij}+l_j^2\U_{i,j-2}}{l_{j}l_{j-1}(l_{j}+l_{j-1})}.
\end{equation*}
\end{itemize}
By combining $(i),(ii),(iii),(iv)$ in \eqref{quad-upw2},  for $i,j=2,\ldots,N-1$, we have 
\begin{eqnarray}
\label{flux-upw2}
&   &  \\
J^{ij}               & = & \epsilon_{ij}\U_{i-2,j}+\eta_{ij}\U_{i-1,j}+\kappa_{ij}\U_{i,j-2}
               +\mu_{ij}\U_{i,j-1}+\Omega_{ij}\U_{ij}+\phi_{ij}\U_{i,j+1}
               +\Psi_{ij}\U_{i,j+2}+\Delta_{ij}\U_{i+1,j}+\Pi_{ij}\U_{i+2,j} \nonumber
\end{eqnarray}
where
\begin{eqnarray*}
&  & \epsilon_{ij}=\frac{h_i^2}{h_ih_{i-1}(h_i+h_{i-1})}\min(p_i,0)~~~~~~~~~\eta_{ij}=-\frac{(h_i+h_{i-1})^2}{h_ih_{i-1}(h_i+h_{i-1})}\min(p_i,0)~~~~~~~~\kappa_{ij}=\frac{l_j^2}{l_jl_{j-1}(l_j+l_{j-1})}\min(q_j,0)\\
&    & \\
&    & \\
&     & ~~~~~~~~~\mu_{ij}=-\frac{(l_j+l_{j-1})^2}{l_jl_{j-1}(l_j+l_{j-1})}\min(q_j,0)\\
&   & \\
&   & \\
&   & \Omega_{ij}=\omega+\frac{(h_i+h_{i-1})^2-h_i^2}{h_{i}h_{i-1}(h_{i}+h_{i-1})}\min(p_i,0)
+\frac{h_{i+1}^2-(h_{i+1}+h_{i+2})^2}{h_{i+1}h_{i+2}(h_{i+1}+h_{i+2})}\max(p_i,0)+\frac{(l_j+l_{j-1})^2-l_j^2}{l_jl_{j-1}(l_j+l_{j-1})}\min(q_j,0)\\
&    & \\
&   & ~~~~~~~+\frac{l_{j+1}^2-(l_{j+1}+l_{j+2})^2}{l_{j+1}l_{j+2}(l_{j+1}+l_{j+2})}\max(q_j,0)\\
&   & \\
&   & \phi_{ij}=\frac{(l_{j+1}+l_{j+2})^2}{l_{j+1}l_{j+2}(l_{j+1}+l_{j+2})}\max(q_j,0)~~~~~~~
\Psi_{ij}=-\frac{l_{j+1}^2}{l_{j+1}l_{j+2}(l_{j+1}+l_{j+2})}\max(q_j,0)~~~~~~~\\
&    &\\
&    &\\
&   &  \Delta_{ij}=\frac{(h_{i+1}+h_{i+2})^2}{h_{i+1}h_{i+2}(h_{i+1}+h_{i+2})}\max(p_i,0)~~~~~~~~~~
\Pi_{ij}=-\frac{h_{i+1}^2}{h_{i+1}h_{i+2}(h_{i+1}+h_{i+2})}\max(p_i,0).
\end{eqnarray*}
For the control volumes near the boundary of the study domain,  two ghost  points  or the first order upwind method can be used.
Finally, we have the following matrix form
\begin{equation}
\label{flux-up2-mat}
J=A_{2up} \U+F_{2up}
\end{equation}
where 
\begin{equation*}
J=\begin{bmatrix}
J^{11}\\
J^{12}\\
\vdots\\
J^{1N}\\
J^{21}\\
J^{22}\\
\vdots\\
J^{2N}\\
\vdots\\
\vdots\\
J^{NN}\\
\end{bmatrix},
F_{2up}=\begin{bmatrix}
F^{11}_{up}\\
F^{12}_{up}\\
\vdots\\
F^{1N}_{up}\\
F^{21}_{up}\\
F^{22}_{up}\\
\vdots\\
F^{2N}_{up}\\
\vdots\\
\vdots\\
F^{NN}_{up}
\end{bmatrix},~~~~~~
\U=\begin{bmatrix}
\U_{11}\\
\U_{12}\\
\vdots\\
\U_{1N}\\
\U_{21}\\
\U_{22}\\
\vdots\\
\U_{2N}\\
\vdots\\
\vdots\\
\U_{NN}\\
\end{bmatrix}
\end{equation*}

and 
\begin{equation*}
A_{2up}=\begin{bmatrix}
H_1 & P_1 &   0_N & 0 & \ldots & \ldots & \ldots &  & 0_N & 0_N\\ 
Q_2 & H_2 & P_2& R_2 & 0_N &   &  &   & & 0_N \\
W_3& Q_3 & H_3 & P_3 &  R_3 & 0_N &   &    &  & \vdots  \\
0_N & W_4 & Q_4 & H_4 & P_4 & R_4 & \ddots\\
0_N    &    0_N &  \ddots & \ddots &  \ddots & \ddots & \ddots & \ddots \\
\vdots &   &  \ddots & \ddots & \ddots & \ddots & \ddots & \ddots                                          & \ddots \\
\vdots &   &  & \ddots & \ddots & \ddots & \ddots & \ddots                                          & \ddots & 0_N \\
 &   &   &  &  \ddots &  W_{N-2} & Q_{N-2}  & H_{N-2} & P_{N-2}  & R_{N-2} \\
\vdots &  &  &  & & \ddots & W_{N-1} & Q_{N-1} & H_{N-1} & P_{i,N-1}\\
 0_N & \dots   &   \ldots  & \ldots  &   & \ldots & 0_N & 0_N  & Q_{N}  & H_{N} \\
\end{bmatrix}
\end{equation*}

where $H_1,H_N$ are tridiagonal matrices, for $i=2,\ldots,N-1, ~H_i$ are penta-diagonal matrices and $P_i,R_i,W_i,Q_i$ are diagonal matrices, and  $F_{up}$ is a vector coming from the boundary conditions.
A structure of the advection matrix using the second order upwind method can be viewed in \figref{upw22}.
As for the first order upwinding, combining the MPFA method \eqref{flux-mpfa-mat} and the second order upwind method \eqref{flux-up2-mat}, we have 
\begin{equation}
\label{mpfa-up2}
\frac{d\U}{d\tau}=A\U+F
\end{equation} 
\begin{eqnarray*}
	A=L^{-1}\Bigg(A_{mp}+A_{2up}+A_L\Bigg)~~~~~~F=L^{-1}\Bigg(F_{mp}+F_{2up}\Bigg),
\end{eqnarray*}
where $A_L$ is a diagonal matrix of size $N^2\times N^2$ coming from the discretisation of \eqref{linearterm}. 
The elements of $A_L$ are $h_il_j\lambda$ for $i,j=1,...,N$ with $\lambda$ given in \eqref{conservation}.
The matrix $L$ is also a diagonal matrix of size $N^2\times N^2$ whose diagonal elements are $L_{ii}=h_il_i$ for $i=1,\ldots,N$

Actually, the PDE \eqref{twoop} is degenerated when the stock price is approaching zero
 $(x\rightarrow 0, y\rightarrow 0)$ which  has an adverse impact on the accuracy of the numerical method. However, to overcome the degeneracy, we are going to apply a fitted finite volume method in the degeneracy region $(x\rightarrow 0, y\rightarrow 0)$. More details about this fitted method is given in the next section.

\subsection{Fitted Multi-Point Flux Approximation }

The fitted Multi-Point Flux Approximation is a combination of the fitted finite volume method ( see \cite{huang2006fitted,huang2009convergence}) and the Multi-Point Flux Approximation method. 
The fitted finite volume helps to deal with the degeneracy of the PDE \eqref{twoop}. 
We approximate simultaneously the diffusion term and the convection term in the degeneracy region by solving a two-points boundary problem. 
In the region where the PDE is not degenerated, we apply the standard Multi-point flux approximation to the diffusion term as described in the previous section.\\
Let us set
\begin{equation}
\label{diff-conv}
k(\U)=\nabla \cdot (\M\nabla \U+ f\U)
\end{equation}
where $\M$ and $f$ are defined in \eqref{conservation}.
Thereby, we have the following decomposition  over a control volume $\C_{ij}$, for $~~i,j=1,...,N$
\begin{eqnarray}
\label{diff-conv-int}
\int_{\C_{ij}}\nabla k(\U)d\C &  =  & \int_{\C_{ij}}\nabla \cdot (M\nabla\U+f\U)d\C \nonumber\\
&    &\nonumber \\
& = & \int_{\partial \C_{ij}}(M\nabla\U+f\U) \cdot\vec n d\partial\C\nonumber\\
&   & \nonumber\\
& = & \int_{(x_{i+\frac{1}{2}},y_{j-\frac{1}{2}})}^{(x_{i+\frac{1}{2}},y_{j+\frac{1}{2}})}\Bigg(m_{11}\frac{\partial \U}{\partial x}+m_{12}\frac{\partial \U}{\partial y}+p\U\Bigg)dy\\
&    & \nonumber\\
&    & -\int_{(x_{i-\frac{1}{2}},y_{j-\frac{1}{2}})}^{(x_{i-\frac{1}{2}},y_{j+\frac{1}{2}})}\Bigg(m_{11}\frac{\partial \U}{\partial x}+m_{12}\frac{\partial \U}{\partial y}+p\U\Bigg)dy\nonumber\\
&   & \nonumber\\
&   & +\int_{(x_{i-\frac{1}{2}},y_{j+\frac{1}{2}})}^{(x_{i+\frac{1}{2}},y_{j+\frac{1}{2}})}\Bigg(m_{21}\frac{\partial \U}{\partial x}+m_{22}\frac{\partial \U}{\partial y}+q\U\Bigg)dx\nonumber\\
&   & \nonumber\\
&   & -\int_{(x_{i-\frac{1}{2}},y_{j-\frac{1}{2}})}^{(x_{i+\frac{1}{2}},y_{j-\frac{1}{2}})}\Bigg(m_{21}\frac{\partial \U}{\partial x}+m_{22}\frac{\partial \U}{\partial y}+q\U\Bigg)dx \nonumber
\end{eqnarray}
with $\vec{n}$ is the outward unit normal vector, $m_{11},m_{12},m_{21},m_{22}$ the coefficients of the matrix $\M$ and $p,q$ coefficients of vector $f$ defined in \eqref{conservation}.\\
 In their work, \cite{huang2006fitted,huang2009convergence} showed how the fitted finite method is used to approximate each of the integral in \eqref{diff-conv-int}.
\subsubsection{Fitted Finite volume method in the degeneracy region}
Following \cite{huang2006fitted}, the fitted finite volume method is used to approximate the flux through the edges which are effectively in the degeneracy region notably the western edge of the control volume $\C_{1,j}$ for $j=1,\ldots,N$ and the southern edge of the control volume $\C_{i,1}$ for $i=1,\ldots,N$ .\\
Thereby, the  flux through the southern edge of the control volume $\C_{i,1}$ 
for $i=1,\ldots,N$ is calculated as follows.\\
\begin{figure}[!h]
\centering
\begin{tikzpicture}[scale=0.6]
\draw[black,thick] (0,0)--(0,8)--(16,8)--(16,0)--(0,0);
\draw[black,thick] (4,0)--(4,8);
\draw[black,thick] (8,0)--(8,8);
\draw[black,thick] (12,0)--(12,8);
\draw[black,thick] (16,0)--(16,8);
\draw[black,thick] (0,4)--(16,4);
\draw[blue,dotted,thick] (6,2)--(6,6)--(10,6)--(10,2)--(6,2);
\draw[green,thick] (6,2)--(10,2);
\draw[red,fill=red] (8,4) circle (0.1);
\node[below,right,red] at (8,4.7){$(x_i,y_1)$};
\draw[green,->] (9,1.5)--(9,2.5);
\end{tikzpicture}
\caption{}
\end{figure}
%
%
%
%
The fitted finite volume method is applied to approximate the integral along the southern edge of control volume $\C_{i,1}$. The idea is to approximate the integral over $[x_{i-\frac{1}{2}};x_{i+\frac{1}{2}}]$ by a constant. We start by applying the mid-quadrature rule as follows:
\begin{equation}
\label{south-appr}
\int_{(x_{i-\frac{1}{2}},y_{\frac{1}{2}})}^{(x_{i+\frac{1}{2}},y_{\frac{1}{2}})}\Bigg(m_{21}\frac{\partial \U}{\partial x}+m_{22}\frac{\partial \U}{\partial y}+q\U\Bigg)dx\approx \Bigg(m_{21}\frac{\partial \U}{\partial x}+m_{22}\frac{\partial \U}{\partial y}+q\U\Bigg)_{\vert_{x_i,y_{\frac{1}{2}}}}\cdot h_i
\end{equation}

Besides we have

\begin{equation}
m_{21}\frac{\partial \U}{\partial x}+m_{22}\frac{\partial \U}{\partial y}+q\U=y\Bigg(ey\frac{\partial \U}{\partial y}+h'\frac{\partial \U}{\partial x}+k\U\Bigg)
\end{equation}
with $e=\frac{1}{2}\sigma_2^2,~~~k=r-\sigma_2^2-\frac{1}{2}\rho\sigma_1\sigma_2$   and $h'=\frac{1}{2}\rho\sigma_1\sigma_2x$.\\
\\
We want to approximate 
\begin{equation*}
f(\U)=ey\frac{\partial \U}{\partial y}+k\U
\end{equation*}
by a linear function over $I_{y_1}=(0,y_{1})$ satisfying the following two-points boundary value problem
\begin{align}
\label{two-bvp}
\left\lbrace \begin{array}{l}
f'(\U)~~~~=\Bigg(ey\frac{\partial \U}{\partial y}+k\U\Bigg)'=K_1\\
\\
\U(x_i,0) =\U_{i,0}~~~~~~~~~\U(x_i,y_1)=\U_{i,1}
\end{array}\right.
\end{align}
By solving this problem we get 
\begin{eqnarray}
\label{sol-bvp}
\U=\U_{i,0}+(\U_{i,1}-\U_{i,0})\frac{y}{y_{1}}
\end{eqnarray}
Thereby, by using \eqref{south-appr}, \eqref{two-bvp}, \eqref{sol-bvp} and  the forward difference for approximating the first partial derivative $\frac{\partial \U}{\partial x}$ we get 
\begin{equation}
\label{south-appr-int}
\int_{(x_{i-\frac{1}{2}},y_{\frac{1}{2}})}^{(x_{i+\frac{1}{2}},y_{\frac{1}{2}})}\Bigg(m_{21}\frac{\partial \U}{\partial x}+m_{22}\frac{\partial \U}{\partial y}+q\U\Bigg)dx  \approx  \frac{1}{2}y_1\Big[\frac{1}{2}h_i(e+k)-h_i'\Big]\U_{i,1}+\frac{1}{2}h_i'y_1\U_{i+1,1}-\frac{1}{4}y_1h_i(e-k)\U_{i,0}
\end{equation}
where
\begin{eqnarray*}
e=\frac{1}{2}\sigma_2^2,~~~~~~~~~k=r-\sigma_2^2-\frac{1}{2}\rho\sigma_1\sigma_2  ~~~~~~~h_i'=\frac{1}{2}\rho\sigma_1\sigma_2x_i~~~~~~h_i=x_{i+\frac{1}{2}}-x_{i-\frac{1}{2}}
\end{eqnarray*}
Similarly, for the western edge of the control volume $\C_{1,j},~~for~~j=1,...,N$, we have
\begin{equation}
\label{west-appr-int}
\int_{(x_{\frac{1}{2}},y_{j-\frac{1}{2}})}^{(x_{\frac{1}{2}},y_{j+\frac{1}{2}})}\Bigg(m_{11}\frac{\partial \U}{\partial x}+m_{12}\frac{\partial \U}{\partial y}+p\U\Bigg)dy \approx \frac{1}{2}x_1\Big[\frac{1}{2}l_j(a+b)-d_j\Big]\U_{1,j}+\frac{1}{2}d_jx_1\U_{1,j+1}-\frac{1}{4}l_jx_1(a-b)\U_{0,j}
\end{equation}
with
\begin{eqnarray*}
a=\frac{1}{2}\sigma_1^2~~~~~~~~~b=r-\sigma_1^2-\frac{1}{2}\rho\sigma_1\sigma_2~~~~~~~~~d_j=\frac{1}{2}\rho\sigma_1\sigma_2y_j~~~~~~~~l_j=y_{j+\frac{1}{2}}-y_{j-\frac{1}{2}}
\end{eqnarray*}
 
\subsubsection{Fitted Multi-Point Flux Approximation }

The fitted Multi-Point Approximation method consists of calculating the flux through the edges which are totally in the degeneracy region using the fitted finite volume method 
as described in the previous paragraph. For the edges which are not totally in the degeneracy region, the flux is approximated using simultaneously the Multi-point 
flux approximation and the upwind methods (first order or second order). In the other hand, the MPFA method and the upwind methods are used to approximate 
respectively the diffusion term and the convection term over the control volumes which are not in the degeneracy region. \\
Considering \eqref{diff-conv-int}, 
in fact, in the control volume $\C_{11}$, the southern and western edges are in the degeneracy region, the northern and the eastern edges are not in the degeneracy region. Thereby, the flux through the southern and western edges  are approximated using the fitted finite volume method, while the 
flux through the eastern and northern edges are approximated using simultaneously of the MPFA method and the upwind method.
This gives
\begin{eqnarray}
\label{appr-intc11}
\int_{{\mathcal{C}}_{11}}\nabla k(\U) & \approx & a_{11}^1\U_{11}+b_{11}^1\U_{21}+c_{11}^1\U_{22}+d_{11}^1\U_{12}
+\omega_{11}^1\U_{02}+\phi_{11}^1\U_{01} \nonumber\\
&   & \nonumber\\
&   & +r_{11}^1\U_{10}+s_{11}^1\U_{20}
\end{eqnarray}
with
\begin{eqnarray*}
&& a_{11}^1=T_{11}^{22}+T_{23}^{21}+T_{31}^{22}+T_{42}^{12}+l_1\max(f_x^2,0)
      +h_1\max(f_y^2,0)-\frac{1}{2}x_1\Big[\frac{1}{2}l_1(a+b)-d_1\Big]\\\
      &   & \\
&   & ~~~~~~~~~-\frac{1}{2}y_1\Big[\frac{1}{2}h_1(e+k)-h_1'\Big]\\
&  & \\
&  & b_{11}^1=T_{12}^{22}
+T_{24}^{21}+T_{32}^{22}+l_1\min(f_x^2,0)-\frac{1}{2}h_1'y_1;~~~~~~~~c^1_{11}=T_{14}^{22}+T_{34}^{22}\\
&   & \\
&   & d_{11}^1=T_{13}^{22}+T_{33}^{22}+T_{44}^{12}+h_1\min(f_y^2,0)-\frac{1}{2}d_1x_1;~~~~~~~~ \omega_{11}^1=T_{43}^{12}~~~~~~~~~~~\\
&   & \\
&   & 
\phi_{11}^1=T_{41}^{12}+\frac{1}{4}l_1x_1(a-b)~~~~~~~~~~~r_{11}^1=T_{21}^{21}+\frac{1}{4}h_1y_1(e-k)~~~~~~~s_{11}^1=T_{22}^{21}
\end{eqnarray*}
Similarly, for the control volume $\C_{1,j}~~~~j=1,\ldots,N$, we have
\begin{eqnarray}
\label{appr-intc1j}
\int_{{\mathcal{C}}_{1,j}}\nabla k(\U)      & \approx &  a_{1,j}^1\U_{1,j}+bb_{1,j}\U_{2,j}+c_{1,j}^1\U_{2,j+1}+d_{1,j}^1\U_{1,j+1}
+\gamma_{1,j}^1\U_{1,j-1}+\lambda_{1,j}^1\U_{2,j-1} \nonumber\\
&   & \nonumber\\
&   & \omega_{1,j}^1\U_{0,j+1}+\phi_{1,j}^1\U_{0,j}
+\Upsilon_{1,j}^1\U_{0,j-1}
\end{eqnarray}
\begin{eqnarray*}
&   & a_{1,j}^1=T_{11}^{2,j+1}+T_{23}^{2,j}+T_{31}^{2,j+1}+T_{42}^{1,j+1}-T_{33}^{2,j}-T_{44}^{1,j}-\frac{1}{2}x_1\Big(\frac{1}{2}l_j(a+b)-d_j\Big)\\
      &  & \\
      &  &~~~~~~~~~+l_j\max(f_x^2,0)+h_1\max(f^{j+1}_y,0)-h_1\min(f^j_y,0)\\
&  &  \\
&   & b_{1,j}^1=T_{12}^{2,j+1}+T_{24}^{2,j}+T_{32}^{2,j+1}-T_{34}^{2,j}+l_j\min(f_x^2,0);~~~~~~~~~~~~c_{1,j}^1=T_{14}^{2,j+1}+T_{34}^{2,j+1};\\
&   & \\
&   &d_{1,j}^1=T_{13}^{2,j+1}+T_{33}^{2,j+1}+T_{44}^{1,j+1}+h_1\min(f_y^{j+1},0)-\frac{1}{2}d_jx_1+\\
&    & \\
&   & \gamma_{1,j}^1=T_{21}^{2,j}-T_{31}^{2,j}-T_{42}^{1,j}-h_1\max(f_y^{j},0);~~~~~~~~~~~~~\lambda_{1,j}^1=T_{22}^{2,j}-T_{32}^{2,j};\\
&   & \\
&   & \omega_{1,j}^1=T_{43}^{1,j+1};~~~~~~~~~~~~~~\phi_{1,j}^1=T_{41}^{1,j+1}-T_{43}^{1,j}+\frac{1}{4}l_jx_1(a-b);~~~~~~~~
\Upsilon_{1,j}^1=-T_{41}^{1,j};
\end{eqnarray*}
For the control $\C_{i,1}~~~~i=2,\ldots,N$, we have:
\begin{eqnarray}
\label{appr-intci1}
\int_{{\mathcal{C}}_{i,1}}\nabla k(\U) & \approx & a_{i,1}^1\U_{i,1}+b_{i,1}^1\U_{i+1,1}+c_{i,1}^1\U_{i+1,2}+d_{i,1}^1\U_{i,2}+e_{i,1}^1\U_{i-1,2}
+\alpha_{i,1}^1\U_{i-1,1}+t_{i,1}^1\U_{i-1,0} \nonumber\\
&   & \nonumber\\
&   & +r_{i,1}^1\U_{i,0}+s_{i,1}^1\U_{i+1,0}
\end{eqnarray}
with 
\begin{eqnarray*}
& &  a_{i,1}^1=T_{11}^{i+1,2}+T_{23}^{i+1,1}+T_{31}^{i+1,2}+T_{42}^{i,2}-T_{12}^{i,2}-T_{24}^{i,1}
-\frac{1}{2}y_1\Big[\frac{1}{2}h_i(e+k)-h_i'\Big]\\
&   &  \\
&    & ~~~~~~~~~~~+l_1\max(f_x^{i+1},0)+h_i\max(f_y^2,0)\Big)-l_1\min(f_x^{i},0)\\
&  & \\
&   & b_{i,1}^1=T_{12}^{i+1,2}+T_{24}^{i+1,1}+T_{32}^{i+1,2}+l_i\min(f_x^{i+1},0)-\frac{1}{2}h_i'y_1;~~~~~~~~~~~~~~~c_{i,1}^1=T_{14}^{i+1,2}+T_{34}^{i+1,2}\\
&   & \\
&   & d_{i,1}^1=T_{13}^{i+1,2}+T_{33}^{i+1,2}+T_{44}^{i,2}-T_{14}^{i,2}
                      +h_i\min(f_y^2,0);~~~~~~~~~~~~~~~~~~~~~~e_{i,1}^1=T_{43}^{i,2}-T_{13}^{i,2}\\
                      &    & \\
                      &    & \alpha_{i,1}^1=T_{41}^{i,2}-T_{11}^{i,2}-T_{23}^{i,1}-l_1\max(f_x^{i},0);~~~~~~~~~~~t_{i,1}^1=-T_{21}^{i,1};\\
&   & \\
&   & r_{i,1}^1=T_{21}^{i+1,1}-T_{22}^{i,1}+\frac{1}{4}y_1h_i(e-k)~~~~~~~~s_{i,1}^1=T_{22}^{i+1,1}
\end{eqnarray*}
As we already mentioned, for  the control volumes which are not in the degeneracy region, we use the multi-Point flux approximation 
to approximate the diffusion term and the upwind methods (first and second order) to approximate the convection  term.
So by combining as before,
we obtain the following ODE
\begin{equation}
\label{fit-mpfa-up1}
\frac{d\U}{d\tau}=A\U+F
\end{equation} 
where
\begin{equation*}
\U=\begin{bmatrix}
\U_{11}\\
\U_{12}\\
\vdots\\
\U_{1N}\\
\U_{21}\\
\U_{22}\\
\vdots\\
\U_{2N}\\
\vdots\\
\vdots\\
\U_{N,1}\\
\U_{N,2}\\
\vdots\\
\U_{NN}
\end{bmatrix}
~~~~A=L^{-1}\Big(Z+A_L\Big)~~
\end{equation*}

with $F$ the vector of boundary conditions, $A_L$ is a diagonal matrix of size $N^2\times N^2$ coming from the discretisation of \eqref{linearterm}.
The elements of $A_L$ are $h_il_j\lambda$ for $ i,j=1,...,N$ with $\lambda$ given in \eqref{conservation}. 
The matrix $L$ is also a diagonal matrix of size $N^2\times N^2$ whose diagonal elements are $h_il_j$ for $i,j=1,\ldots,N$
 and
\begin{equation*}
Z=\begin{bmatrix}
D_1 & K_1 & 0_N & \ldots & \ldots & \ldots & \ldots & 0_N\\ 
L_2 & D_2 & K_2 & \ddots &  &  &  & \vdots  \\ 
0_N & L_3 & D_3 & K_3 & \ddots &  & & \vdots\\ 
\vdots & \ddots & L_4 & D_4 & K_4 & \ddots & & \vdots\\ 
\vdots &  & \ddots & \ddots & \ddots & \ddots & \ddots & \vdots \\ 
\vdots &  &  & \ddots& \ddots & \ddots & \ddots & 0_N&  \\ 
\vdots &  &  &  & \ddots & L_{N-1} & D_{N-1} & K_{N-1} \\ 
0_N & \ldots & \ldots & \ldots & \ldots & 0_N & L_N & D_N
\end{bmatrix} 
\end{equation*}
The fitted matrix $Z$ uses the first order upwind method.
The matrices $D_i,K_i,L_i$ are tri-diagonal matrices  defined as follows. For $i=1,N$
\begin{eqnarray*}
k=1,\ldots,N~~(D_i)_{kk}=a_{1,k}^1~~~~~~~~~k=1,\ldots,N-1~~(D_i)_{k,k+1}=d_{1,k}^1,~~~~~~~~~~~~
k=2,\ldots,N~~(D_i)_{k,k-1}=\gamma_{1,k}^1\\
&  & \\
k=1,\ldots,N~~(K_1)_{kk}=b_{1,k}^1~~~~~~~~~k=1,\ldots,N-1~~(K_1)_{k,k+1}=c_{1,k}^1,~~~~~~~~~~~~
k=2,\ldots,N~~(K_1)_{k,k-1}=\lambda_{1,k}^1\\
&   & \\
(L_N)_{11}=\alpha_{N,1}^1~~~(L_N)_{12}=e_{N,1}^1~~~~~~~~~~~~~~~~~~~~~~~~~~~
~~~~~~~~~~~~~~~~~~~~~~~~~~~~~~~~~~\\
&   & \\
k=2,\ldots,N~~(L_N)_{kk}=\alpha_{N,k}+\epsilon_{N,k}~~~~~~~~~k=1,\ldots,N-1~~(L_N)_{k,k+1}=e_{N,k},~~~~~~~~~k=2,\ldots,N~~(L_N)_{k,k-1}=\beta_{N,k}
\end{eqnarray*}

For $i=2,\ldots,N-1$

\begin{eqnarray*}
&   & \\
 &  & ~(D_i)_{11}=a_{i,1}^1~;~(D_i)_{12}=d_{i,1}^1;~~~~~(K_i)_{11}=b_{i,1}^1~
 ;~(K_i)_{12}=c_{i,1}^1~~~~(L_i)_{11}=\alpha_{i,1}~;~(L_i)_{12}=e_{i,1}^1\\
&   & \\
&   & 
k=2,\ldots,N~~~~~(D_i)_{kk}=a_{i,k}+\Omega_{i,k};~~~~~~~
(K_i)_{kk}=b_{i,k}+\psi_{i,k};~~~~~~~~
(L_i)_{kk}=\alpha_{i,k}+\epsilon_{i,k}\\
&   & \\
&   & k=2,\ldots,N-1~~~~~(D_i)_{k,k+1}=d_{i,k}+\phi_{i,k};~~~~~~~
(K_i)_{k,k+1}=c_{i,k};~~~~~~~~
(L_i)_{k,k+1}=e_{i,k}\\
&   &  \\
&   & k=2,\ldots,N~~~~~(D_i)_{k,k-1}=\gamma_{i,k}+\mu_{i,k};~~~~~~~
(K_i)_{k,k-1}=\lambda_{i,k};~~~~~~~~
(L_i)_{k,k-1}=\beta_{i,k}
\end{eqnarray*}
where all the elements  $a_{i,j}^1,b_{i,j}^1,c_{i,j}^1,d_{i,j}^1,e_{i,j}^1,\gamma_{i,j}^1,\lambda_{i,j}^1$ are defined in \eqref{appr-intc11},\eqref{appr-intc1j},\eqref{appr-intci1} and the others elements are defined in \eqref{flux-mpfa} and \eqref{flux-up1}.

Similarly, combining the fitted finite volume method, the MPFA and the second  order upwind method 
 we have
\begin{equation}
\label{fit-mpfa-up2}
\frac{d\U}{d\tau}=A\U+F
\end{equation} 
where
\begin{equation*}
\U=\begin{bmatrix}
\U_{11}\\
\U_{12}\\
\vdots\\
\U_{1N}\\
\U_{21}\\
\U_{22}\\
\vdots\\
\U_{2N}\\
\vdots\\
\vdots\\
\U_{N,1}\\
\U_{N,2}\\
\vdots\\
\U_{NN}
\end{bmatrix}
~~~~A=L^{-1}\Big(Y+A_L\Big)~~
\end{equation*}

with $G$ the vector of boundary conditions, $A_L$ is a diagonal matrix of size $N^2\times N^2$ coming from the discretisation of \eqref{linearterm}. The elements of $A_L$ are $h_il_j\lambda$ for $i,j=1,...,N$ with $\lambda$ given in \eqref{conservation}. The matrix L is also a diagonal matrix of size $N^2\times N^2$ whose elements are $h_il_j$ for $i,j=1,\ldots,N$
and

\begin{equation*}
Y=\begin{bmatrix}
H_1 & P_1 &   0_N & 0 & \ldots & \ldots & \ldots &  & 0_N & 0_N\\ 
Q_2 & H_2 & P_2& R_2 & 0_N &   &  &   & & 0_N \\
W_3& Q_3 & H_3 & P_3 &  R_3 & 0_N &   &    &  & \vdots  \\
0_N & W_4 & Q_4 & H_4 & P_4 & R_4 & \ddots\\
0_N    &    0_N &  \ddots & \ddots &  \ddots & \ddots & \ddots & \ddots \\
\vdots &   &  \ddots & \ddots & \ddots & \ddots & \ddots & \ddots                                          & \ddots \\
\vdots &   &  & \ddots & \ddots & \ddots & \ddots & \ddots                                          & \ddots & 0_N \\
 &   &   &  &  \ddots &  W_{N-2} & Q_{N-2}  & H_{N-2} & P_{N-2}  & R_{N-2} \\
\vdots &  &  &  & & \ddots & W_{N-1} & Q_{N-1} & H_{N-1} & P_{i,N-1}\\
 0_N & \dots   &   \ldots  & \ldots  &   & \ldots & 0_N & 0_N  & Q_{N}  & H_{N} \\
\end{bmatrix}
\end{equation*}

%
The elements of matrix Y are matrices. Indeed $0_N$ is a zeros matrix of size $N\times N$. The matrices $H_i,P_i,Q$ are tri-diagonal matrices and $W_i,R_i$ are diagonal matrices  defined as follows:

\begin{eqnarray*}
&  & (H_1)_{11}=a_{11}^1~;~(H_1)_{12}=d_{11}^1~~~~~~(P_1)_{11}=b_{11}^1~;~
(P_1)_{12}=c_{11}^1\\
&   & \\
&   & k=2,\ldots, N~~(H_1)_{kk}=a_{1,k}^1;~~~~~k=2,\ldots, N-1~~(H_1)_{k,k+1}=d_{1,k}^1;~~~~~~k=2,\ldots ,N
~~(H_1)_{k,k-1}=\gamma_{1,k}^1\\
&   & \\
&   & k=2,\ldots, N~~(P_1)_{kk}=b_{1,k}^1;~~~~~k=2,\ldots, N-1~~(P_1)_{k,k+1}=c_{1,k}^1;~~~~~~k=2,\ldots ,N
~~(P_1)_{k,k-1}=\lambda_{1,k}^1
\end{eqnarray*}

For $i=2,\ldots,N-1$

\begin{eqnarray*}
 &  & ~(H_i)_{11}=a_{i,1}^1~;~(H_i)_{12}=d_{i,1}^1;~~~~~(P_i)_{11}=b_{i,1}^1+\Delta_{i,1}~
 ;~(P_i)_{12}=c_{i,1}^1~~~~(Q_i)_{11}=\alpha_{i,1}+\eta_{i,1}~;~(Q_i)_{12}=e_{i,1}^1\\
&   & \\
&   & 
k=2,\ldots, N,~~~(H_i)_{kk}=a_{i,k}+\Omega_{i,k};~~~~~~~~
(P_i)_{kk}=b_{i,k}+\Delta_{i,k};~~~~~~~~~
(Q_i)_{kk}=\alpha_{i,k}+\eta_{i,k}~~~~~~
\\
&   & \\
&   & k=2,\ldots,N-1, ~~~~~ (H_i)_{k,k+1}=d_{i,k}+\phi_{i,k};~~~~~~~
(P_i)_{k,k+1}=c_{i,k};~~~~~~~~
(Q_i)_{k,k+1}=e_{i,k}\\
&   &  \\
&   & k=2,\ldots,N,~~~~~(H_i)_{k,k-1}=\lambda_{i,k}+\mu_{i,k};~~~~~~~
(P_i)_{k,k-1}=\lambda_{i,k};~~~~~~~~
(Q_i)_{k,k-1}=\beta_{i,k}\\
&  & \\
&  & k=2,\ldots, N-2, ~~~~~(H_i)_{k,k+2}=\Psi_{i,k};~~~~~~~~~~~~~~~~~k=3,\ldots,N~~~(H_i)_{k,k-2}=\kappa_{i,k}\\
&   & \\
\end{eqnarray*}
and 
\begin{eqnarray*}
&   &~~ (R_i)_{kk}=\Pi_{ik},~~i=2,\ldots, N-2,\,\,k=2,\ldots, N-1\\
&   &~ (W_i)_{kk}=\epsilon_{ik},~~~i=3,\ldots, N-1, ~~=2,\ldots, N-1,
\end{eqnarray*}

where all the elements  $a_{i,j}^1,b_{i,j}^1,c_{i,j}^1,d_{i,j}^1,e_{i,j}^1,\gamma_{i,j}^1,\lambda_{i,j}^1$ are defined \eqref{appr-intc11},\eqref{appr-intc1j},\eqref{appr-intci1}, and the others elements are defined in \eqref{flux-mpfa} and \eqref{flux-upw2}.

\section{Time discretization}

Let us consider the ODE stemming from the spatial dicretization and given by \eqref{mpfa-up1},\eqref{mpfa-up2},\eqref{fit-mpfa-up1} and \eqref{fit-mpfa-up2}

\begin{equation*}
\frac{d\U}{d\tau}=A\U+F
\end{equation*}

Using the $\theta$-method for the time discretization, we have

\begin{eqnarray}
\frac{\U^{n+1}-\U^n}{\Delta\tau}=\theta\Big(A\U^{n+1}+F^{n+1}\Big)+(1-\theta)\Big(A\U^n+F^n\Big)
\end{eqnarray}

Hence

\begin{equation}
\U^{n+1}=\Big(I-\theta\Delta\tau A\Big)^{-1}\Bigg[\Big(I+(1-\theta)\Delta\tau A\Big)\U^n+\theta \Delta \tau F^{n+1}+(1-\theta)\Delta\tau F^n\Bigg]
\end{equation}

with

\begin{eqnarray*}
&  & \U^n=\begin{bmatrix}
\U_{11}(\tau_n)~~
\U_{12}(\tau_n)~~
\ldots~~
\U_{1N}(\tau_n)~~
\U_{21}(\tau_n)~~
\U_{22}(\tau_m)~~
\ldots~~
\U_{2N}(\tau_n)~~
\ldots~~
\U_{N,1}(\tau_n)~~
\U_{N,2}(\tau_n)~
\ldots\ldots~
\U_{NN}(\tau_n)
\end{bmatrix}^T\\
&  & \\
&  &
 F^n=F(\tau_n),\,\,\;\;\tau_n=n\Delta \tau.
\end{eqnarray*}
\section{Numerical experiments}

In this section, we perform some numerical simulations and compare different numerical schemes developed in this work. More precisely, we compare
the novel fitted  MPFA method combined to the upwind methods, first method (fitted  MPFA-$1^{st}$ upw)  and second order (fitted MPFA-$2^{nd}$ upw),
with  the fitted  finite volume method by  \cite{huang2006fitted} (fitted FV) and   
 the standard MPFA method combined to the upwind methods, first (MPFA-$1^{st}$ upw) and second order (MPFA-$2^{nd}$ upw). The analytical solution of the PDE (\ref{twoop}) is well known (see \cite{haug2007complete} ) and given as

\begin{eqnarray}
\label{analsol}
C(x,y,K,T) & = &  xe^{-rT}M(y_1,d;\rho_1)+ye^{-rT}M(y_2,-d+\sigma\sqrt{T},\rho_2) \nonumber\\
    &   &   \\
    &    & -Ke^{-rT}\times\left(1-M(-y_1+\sigma_1\sqrt{T},-y_2+\sigma_2\sqrt{T},\rho)\right) \nonumber\
\end{eqnarray} 
where
\begin{eqnarray*}
&  & d=\frac{\ln(x/y)+(b_1-b_2+\sigma_1^2/2)T}{\sigma\sqrt{T}},\\
&   & \\
    &   & y_1  =  \frac{\ln(x /K)+(b_1+\sigma_1^2/2)T}{\sigma_1\sqrt{T}},~~~~~~y_2=\frac{\ln(y/K)+(b_1+\sigma_2^2/2)T}{\sigma_2\sqrt{T}},\\
    &   &   \\
    &   & \sigma=\sqrt{\sigma_1^2+\sigma_2^2-2\rho\sigma_1\sigma_2},~~~~~\rho_1=\frac{\sigma_1-\rho\sigma_2}{\sigma}~~~~~~~\rho_2=\frac{\sigma_2-\rho\sigma_1}{\sigma},
\end{eqnarray*}
and
\begin{equation*}
M(a,b,\rho)=\frac{1}{2\pi\sqrt{1-\rho^2}}\int_{-\infty}^a\int_{-\infty}^b \exp\left(-\frac{u^2-2\rho uv+v^2}{2(1-\rho^2)}\right)dudv.
\end{equation*}
Note that in all  our numerical schemes,  the Dirichlet Boundary condition is used with the value equal to the analytical solution.

\begin{figure}[hbtp]
\centering
\includegraphics[scale=0.4]{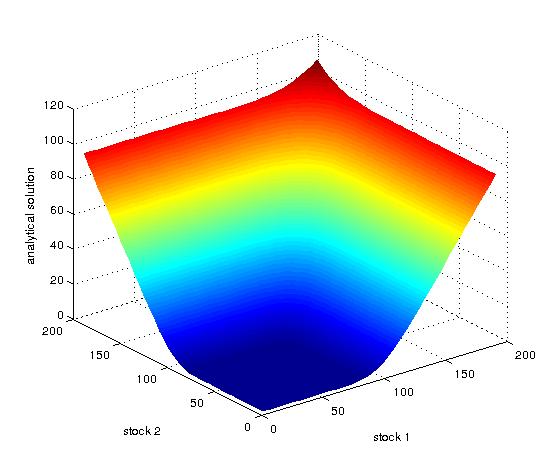}
\caption{ Analytical solution for option price at final time $T$.  The computational domain of the problem is $\Omega=[0;300]\times [0;300]\times[0;T]$ with $T=1/12$, $K=100$, the volatilities $\sigma_1=\sigma_2=0.3$.
The correlation coefficient  is $\rho=0.5$, the risk free interest $r=0.03$ and $\Delta  \tau=1/100$.}
\label{fig1}
\end{figure}
\begin{figure}[hbtp]
\centering
\subfloat[MPFA-upwind 1st order]{
\includegraphics[scale=0.3]{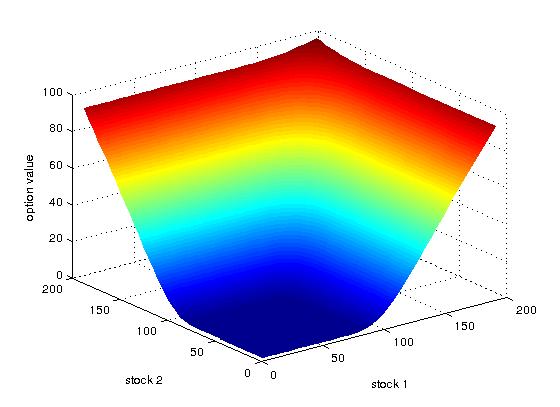}}\qquad
\subfloat[MPFA-upwind 2nd order]{
\includegraphics[scale=0.3]{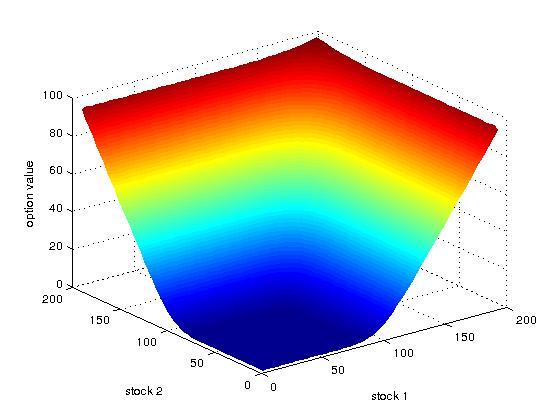}}
\caption{ Option price for MPFA-upwind methods at final time $T$. The computational domain of the problem is $\Omega=[0;300]\times [0;300]\times[0;T]$ with $T=1/12$, $K=100$, the volatilities $\sigma_1=\sigma_2=0.3$.
The correlation coefficient  is $\rho=0.5$, the risk free interest $r=0.03$ and $\Delta  \tau=1/100$.}
\label{fig2}
\end{figure}
\begin{figure}[hbtp]
\centering
\subfloat[fitted MPFA-upwind 1st order]{
\includegraphics[scale=0.3]{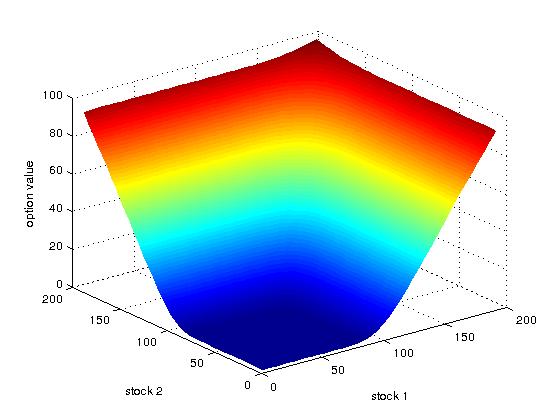}}\qquad
\subfloat[fitted MPFA-upwind 2nd order]{
\includegraphics[scale=0.3]{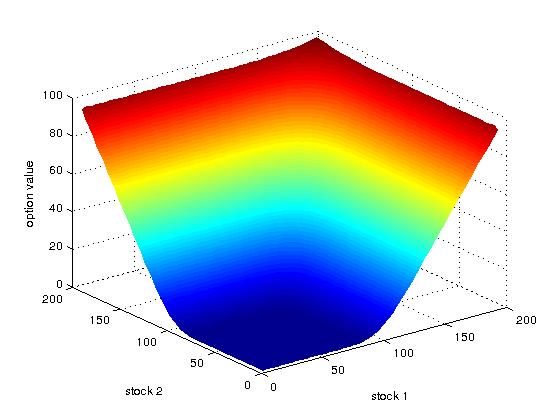}}
\caption{Option price for fitted MPFA-upwind methods at final time $T$.The computational domain of the problem is $\Omega=[0;300]\times [0;300]\times[0;T]$ with $T=1/12$, $K=100$, the volatilities $\sigma_1=\sigma_2=0.3$.
The correlation coefficient  is $\rho=0.5$, the risk free interest $r=0.03$ and  $\Delta  \tau=1/100$. }
\label{fig3}
\end{figure}

The graphs of option price with  different methods are given  in \figref{fig1},\figref{fig2} and \figref{fig3} 

\newpage

In this paragraph, we consider the four numerical methods illustrated in the previous sections and the fitted finite volume method \cite{huang2006fitted}. 
We evaluate the error of these numerical method with respect to the analytical solution \eqref{analsol}. The $L^2$-norm is used to compute the error as follows:

\begin{equation}
\label{error-l2}
err=\frac{\sqrt{\sum_{i,j=1}^N meas(\C_{ij}) \big(\U_{ij}-U_{ij}^{ana}\big)^2}}{\sqrt{\sum_{i,j=1}^n meas(\C_{ij}) \big(U_{ij}^{ana}\big)^2}}
\end{equation}

where $\U$ is the numerical solution, $U^{ana}$ the analytical solution and $meas(\C_{i,j})$ is the measure of the control volume $\C_{ij}$.
This gives the following table: 

\begin{table}[!h]
\centering
\begin{tabular}{|c||c|c|c|c|c|}
\hline
\backslashbox{Nb of grid pts}{Num method}& Fitted fin vol & MPFA-$1^{st}$
upw
  & MPFA-$2^{nd}$ upw &  fitted MPFA-$1^{st}$ upw & fitted  MPFA -$2^{nd}$ upw\\
  \hline 
  $50\times 50$ & 0.0134&  0.0060& 0.0059
 &  0.0060  & 0.0060 \\
  \hline 
$70\times 70$ & 0.0133 &  0.0044 & 0.0044  &  0.0044 & 0.0044 \\
\hline
$85\times 85$ & 0.0132 & 0.0037     &    0.0037   & 0.0037    &  0.0037\\
\hline
$100\times 100$ & 0.0132 & 0.0032   & 0.0032& 0.0032 & 0.0032 \\ 
\hline 
$150\times 150$ & 0.0131 & 0.0024  & 0.0023& 0.0023 & 0.0023 \\ 
\hline
\end{tabular}
\caption{Table of errors.  The computational domain of the problem is $\Omega=[0;300]\times [0;300]\times[0;T]$ with $T=1/6$, $K=100$, the volatilities $\sigma_1=\sigma_2=0.3$.
The correlation coefficient  is $\rho=0.5$, the risk free interest $r=0.1$ and $\Delta \tau=1/100$.}
\label{errorss1}
\end{table}
\begin{table}[!h]
\centering
\begin{tabular}{|c||c|c|c|c|c|}
\hline
\backslashbox{Nb of grid pts}{Num method}& Fitted fin vol & MPFA-$1^{st}$
upw
  & MPFA-$2^{nd}$ upw &  fitted MPFA-$1^{st}$ upw & fitted  MPFA -$2^{nd}$ upw\\
  \hline 
  $50\times 50$ & 0.0134&  0.0060& 0.0059
 &  0.0060  & 0.0060 \\
  \hline
$100\times 100$ & 0.0104 & 0.0064   & 0.0063& 0.0064 & 0.0063 \\ 
\hline 
$150\times 150$ & 0.0131 & 0.0056 & 0.0055&  0.0056 & 0.0055 \\ 
\hline
\end{tabular}
\caption{Table of errors. The computational domain of the problem is $\Omega=[0;300]\times [0;300]\times[0;T]$ with $T=1/6$, $K=100$, the volatilities $\sigma_1=\sigma_2=0.3$.
The correlation coefficient  is $\rho=0.5$, the risk free interest $r=0.08$ and $\Delta  \tau=1/100$.}
\label{errorss2}
\end{table}
\begin{table}[!h]
\centering
\begin{tabular}{|c||c|c|c|c|c|}
\hline
\backslashbox{Nb of grid pts}{Num method}& Fitted fin vol & MPFA-$1^{st}$
upw
  & MPFA-$2^{nd}$ upw &  fitted MPFA-$1^{st}$ upw & fitted  MPFA -$2^{nd}$ upw\\
  \hline
$100\times 100$ & 0.0152 & 0.0239   & 0.0235& 0.0240 & 0.0229  \\ 
\hline 
$150\times 150$ &  0.0151 & 0.0231 & 0.0228&  0.0232 & 0.0229 \\ 
\hline
\end{tabular}
\caption{Table of errors. The computational domain of the problem is $\Omega=[0;300]\times [0;300]\times[0;T]$ with $T=1/6$, $K=100$, the volatilities $\sigma_1=\sigma_2=0.3$.
The correlation coefficient  is $\rho=0.5$ , the risk free interest $r=0$ and $\Delta  \tau=1/100$. }
\label{errorss3}
\end{table}

\begin{table}[!h]
\centering
\begin{tabular}{|c||c|c|c|c|c|}
\hline
\backslashbox{Nb of grid pts}{Num method}& Fitted fin vol & MPFA-$1^{st}$
upw
  & MPFA-$2^{nd}$ upw &  fitted MPFA-$1^{st}$ upw & fitted  MPFA -$2^{nd}$ upw\\
  \hline
$50 \times 50$ & 0.1208 & 0.0631   & 0.0669& 0.0623 & 0.0659 \\ 
\hline 
$100\times 100$ &  0.1203 & 0.0572 & 0.0648&  0.0559 & 0.0629 \\ 
\hline
\end{tabular}
\caption{Table of errors. The computational domain of the problem is $\Omega=[0;4]\times [0;4]\times[0;T]$ with $T=2$, $K=1$, the volatilities $\sigma_1=\sigma_2=1$.
The correlation coefficient  is $\rho=0.3$, the risk free interest $r=0.5$ and $\Delta  \tau=1/100$.}
\label{errorss3}
\end{table}

\begin{table}[!h]
\centering
\begin{tabular}{|c||c|c|c|c|c|}
\hline
\backslashbox{Nb of grid pts}{Num method}& Fitted fin vol & MPFA-$1^{st}$
upw
  & MPFA-$2^{nd}$ upw &  fitted MPFA-$1^{st}$ upw & fitted  MPFA -$2^{nd}$ upw\\
  \hline
$50 \times 50$ & 0.1196 & 0.0562  & 0.0643& 0.0555 & 0.0624 \\ 
\hline 
$100\times 100$ &  0.1201 & 0.0626 & 0.0664&  0.0618 & 0.0654 \\ 
\hline
\end{tabular}
\caption{Table of errors. The computational domain of the problem is $\Omega=[0;4]\times [0;4]\times[0;T]$ with $T=2$, $K=1$, the volatilities $\sigma_1=\sigma_2=1$.
The correlation coefficient  is $\rho=0.3$, the risk free interest $r=0.5$ and $\Delta  \tau=1/10$.}
\label{errorss4}
\end{table}

As we can observe in \tabref{errorss1}-\tabref{errorss4}, the errors from our fitted MPFA  and  MPFA  methods  are smaller compared to those of fitted finite volume in \cite{huang2006fitted}.
We can also note that when $r$ become smaller,  the gaps between the errors of the fitted finite volume in \cite{huang2006fitted} and  our fitted MPFA  and  MPFA  methods  reduce.

\section{Conclusion}

In this paper, we have presented the Multi-Point Flux Approximation (MPFA) to approximate the diffusion term of Black-Scholes
Partial Differential Equation in its divergence form. The MPFA method coupled with the upwind methods (first and second order)
have been used to solve numerically the Black-Scholes PDE. To handle the degeneracy of Black Scholes PDE, 
we have proposed a novel method based on a combination of the MPFA  method and fitted finite volume by \cite{huang2006fitted}.
We have  performed some numerical simulations, which show that our fitted MPFA method coupled with first or second order upwinding methods are more accurate than the fitted finite volume method by  \cite{huang2006fitted}.
Rigorous convergence proof of the fitted MPFA will be our nearest future  work.

\newpage 

\section*{Acknowledgement}

This work was supported by the Robert Bosch Stiftung through the AIMS ARETE Chair programme (Grant No 11.5.8040.0033.0).



\begin{thebibliography}{}
\bibitem[Aavatsmark(2002)]{aavatsmark2002introduction}
Aavatsmark, I.(2002).
\newblock An introduction to multipoint flux approximations for quadrilateral
  grids.
\newblock \emph{Computational Geosciences}, 6\penalty0 (3-4):\penalty0
  405--432.

\bibitem[Aavatsmark(2007)]{aavatsmark2007multipoint}
Aavatsmark, I. (2007)
\newblock Multipoint flux approximation methods for quadrilateral grids.
\newblock In \emph{9th International forum on reservoir simulation, Abu Dhabi},
  pages 9--13.

\bibitem[Angermann and Wang(2007)]{angermann2007convergence}
Angermann, L. \& Wang, S.(2007).
\newblock Convergence of a fitted finite volume method for the penalized
  black--scholes equation governing european and american option pricing.
\newblock \emph{Numerische Mathematik}, 106\penalty0 (1):\penalty0 1--40.

\bibitem[Bates(1996)]{bates1996jumps}
Bates, D. S.(1996).
\newblock Jumps and stochastic volatility: Exchange rate processes implicit in
  deutsche mark options.
\newblock \emph{The Review of Financial Studies}, 9\penalty0 (1):\penalty0
  69--107.

\bibitem[Duffy(2013)]{duffy2013finite}
 Duffy,  D. J. (2013)
\newblock \emph{Finite Difference methods in financial engineering: A Partial
  Differential Equation approach}.
\newblock John Wiley \& Sons.

\bibitem[Haug(2007)]{haug2007complete}
Haug., E. G. (2007).
\newblock \emph{The complete guide to option pricing formulas}, volume~2.
\newblock McGraw-Hill New York.

\bibitem[Heston(1993)]{heston1993closed}
Heston, S. L. (1993).
\newblock A closed-form solution for options with stochastic volatility with
  applications to bond and currency options.
\newblock \emph{The review of financial studies}, 6\penalty0 (2):\penalty0
  327--343.

\bibitem[Huang et~al.(2006)Huang, Hung, and Wang]{huang2006fitted}
Huang, C-S., Hung,C-H., \& Wang, S.(2006).
\newblock A fitted finite volume method for the valuation of options on assets
  with stochastic volatilities.
\newblock \emph{Computing}, 77\penalty0 (3):\penalty0 297--320.

\bibitem[Huang et~al.(2009)Huang, Hung, and Wang]{huang2009convergence}
Huang, C-S., Hung,C-H., \& Wang, S.(2009).
\newblock On convergence of a fitted finite-volume method for the valuation of
  options on assets with stochastic volatilities.
\newblock \emph{IMA journal of numerical analysis}, 30\penalty0 (4):\penalty0
  1101--1120.

\bibitem[Hull(2003)]{hull2003options}
 Hull, J. C. (2003)
\newblock Options, futures and others.
\newblock \emph{Derivative (Fifth Edition), Prentice Hall}.

\bibitem[Kwok(2008)]{kwok2008mathematical}
 Kwok, Y.-K.(2008).
\newblock \emph{Mathematical models of financial derivatives}.
\newblock Springer.

\bibitem[LeVeque(2004)]{leveque2004finite}
LeVeque, R. J. (2004).
\newblock Finite volume methods for hyperbolic problems.
\newblock \emph{Cambridge Texts in Applied Mathematics}, 39\penalty0
  (1):\penalty0 88--89.

\bibitem[Lie et~al.(2012)Lie, Krogstad, Ligaarden, Natvig, Nilsen, and
  Skaflestad]{lie2012open}
Lie K.-A., Krogstad, S., Ligaarden  I., S. ,
  Natvig ,J. R., Nilsen, H. M., \& B{\aa}rd Skaflestad, B.(2012).
\newblock Open-source matlab implementation of consistent discretisations on
  complex grids.
\newblock \emph{Computational Geosciences}, 16\penalty0 (2):\penalty0 297--322.

\bibitem[Persson and Sydow(2007)]{persson2007pricing}
Persson,   J., \&  Sydow, L. V. (2007).
\newblock Pricing European multi-asset options using a space-time adaptive
  {FD}-method.
\newblock \emph{Computing and Visualization in Science}, 10\penalty0
  (4):\penalty0 173--183.

\bibitem[Sandve et~al.(2012)Sandve, Berre, and Nordbotten]{sandve2012efficient}
Sandve, T. H., Berre, I., \& Nordbotten, J. M.(2012)
\newblock An efficient multi-point flux approximation method for discrete
  fracture--matrix simulations.
\newblock \emph{Journal of Computational Physics}, 231\penalty0 (9):\penalty0
  3784--3800.

\bibitem[Stephansen(2012)]{stephansen2012convergence}
Stephansen, A. F.(2012).
\newblock Convergence of the multipoint flux approximation l-method on general
  grids.
\newblock \emph{SIAM Journal on Numerical Analysis}, 50\penalty0 (6):\penalty0
  3163--3187.

\bibitem[Tambue(2016)]{tambue2016exponential}
Tambue, A. (2016).
\newblock An exponential integrator for finite volume discretization of a
  reaction--advection--diffusion equation.
\newblock \emph{Computers \& Mathematics with Applications}, 71\penalty0
  (9):\penalty0 1875--1897.

\bibitem[Wang(2004)]{wang2004novel}
Wang, S.(2004).
\newblock A novel fitted finite volume method for the black--scholes equation
  governing option pricing.
\newblock \emph{IMA Journal of Numerical Analysis}, 24\penalty0 (4):\penalty0
  699--720.

\bibitem[Wilmott(2005)]{wilmott2005best}
Wilmott, P. (2005)
\newblock \emph{The Best of Wilmott 1: Incorporating the Quantitative Finance
  Review}.
\newblock John Wiley \& Sons.

\bibitem[Wilmott et~al.(1993)Wilmott, Dewynne, and Howison]{wilmott1993option}
Wilmott, P.,  Dewynne, J., \&  Howison , S. (1993).
\newblock \emph{Option pricing: mathematical models and computation}.
\newblock Oxford financial press.
\end{thebibliography}
\end{document}